\documentclass[a4paper,reqno]{amsart}%
\usepackage{amsfonts}
\usepackage{amsmath}
\usepackage{amssymb}
\usepackage{amscd}
\usepackage{graphicx}%
\input xy
\xyoption{all}
\newcounter{numbStep}
\newtheorem{theorem}{Theorem}
\newtheorem*{thmForTheta}{Theorem 6}
\newtheorem{step}[numbStep]{Step}
\newtheorem*{main theorem}{Main Theorem}

\newtheorem{claim}[theorem]{Claim}

\newtheorem{corollary}[theorem]{Corollary}

\newtheorem{lemma}[theorem]{Lemma}

\newtheorem{proposition}[theorem]{Proposition}

\begin{document}
\title{The simplicial volume of closed manifolds covered by $\mathbb{H}^{2}%
\times\mathbb{H}^{2}$}
\author{Michelle Bucher-Karlsson}
\address{Math. Dept, KTH, SE-100 44 Stockholm}
\email{mickar@math.kth.se}
\urladdr{http://www.math.kth.se/\symbol{126}mickar/}
\thanks{Supported by the Swiss National Science Foundation, Grant number PBEZ2-106962.}
\date{}

\begin{abstract}
We compute the value of the simplicial volume for closed, oriented Riemannian
manifolds covered by $\mathbb{H}^{2}\times\mathbb{H}^{2}$ explicitly, thus in
particular for products of closed hyperbolic surfaces. This gives the first
exact value of a nonvanishing simplicial volume for a manifold not admitting a
hyperbolic structure.

\end{abstract}
\maketitle

\section{Introduction}

\setcounter{numbStep}{-1} Our main result is the computation of the Gromov
norm, that is the sup norm, of the Riemannian volume form on the product
$\mathbb{H}^{2}\times\mathbb{H}^{2}$, where $\mathbb{H}^{2}$ denotes the
hyperbolic plane of constant curvature $-1$:

\begin{main theorem}
Let $\omega_{\mathbb{H}^{2}\times\mathbb{H}^{2}}\in H_{c}^{4}(\mathrm{PSL}%
_{2}\mathbb{R}\times\mathrm{PSL}_{2}\mathbb{R},\mathbb{R})$ be the image,
under the Van Est isomorphism, of the Riemannian volume form on $\mathbb{H}%
^{2}\times\mathbb{H}^{2}$. Then%
\[
\left\Vert \omega_{\mathbb{H}^{2}\times\mathbb{H}^{2}}\right\Vert _{\infty
}=\frac{2}{3}\pi^{2}.
\]

\end{main theorem}

Recall that the simplicial volume $\left\Vert M\right\Vert $ of a closed,
oriented manifold $M$ is a topological invariant introduced by Gromov in
\cite{Gr82} and is defined as%
\[
\left\Vert M\right\Vert =\mathrm{inf}\left\{  \left.
\Sigma|a_{\sigma}|\text{ }\right\vert \text{ }\Sigma
a_{\sigma}\sigma\text{ represents the real fundamental class
}[M]\right\} .
\]
As an immediate consequence of our main theorem, we obtain in
Theorem \ref{Thm: simpl vol mfld cov by HxH} below the explicit
proportionality constant relating simplicial volume and volume of
closed, oriented, Riemannian manifolds covered by
$\mathbb{H}^{2}\times\mathbb{H}^{2}$. Indeed, we prove in
\cite[Theorem 2]{Bu06} that the proportionality constant of the
proportionality principle for closed, oriented locally symmetric
space of noncompact type $M^{n}=\Gamma\backslash G/K$ is precisely
the Gromov norm of the volume form in $H_{c}^{n}(G,\mathbb{R})$.
Explicit finite values of the proportionality constant were up to
now only known for hyperbolic manifolds (\cite{Gr82},\cite{Th78}):
It is in this case equal to the maximum volume of ideal geodesic
simplices in $\mathbb{H}^{n}$, a constant which has been computed
explicitly up to dimension $n=6$ only.

\begin{theorem}
\label{Thm: simpl vol mfld cov by HxH}Let $M$ be a closed, oriented Riemannian
manifold whose universal cover $\widetilde{M}$ is isometric to $\mathbb{H}%
^{2}\times\mathbb{H}^{2}$. Then
\[
\left\Vert M\right\Vert =\frac{3}{2\pi^{2}}\mathrm{Vol}(M)\text{.}%
\]

\end{theorem}

In view of Hirzebruch's proportionality principle \cite{Hi56}
relating the volume and the Euler characteristic $\chi$ of locally
symmetric spaces of noncompact type, the conclusion of Theorem
\ref{Thm: simpl vol mfld cov by HxH} can be rewritten as
\[
\left\Vert M\right\Vert =6\chi(M)\text{.}%
\]
In particular, using the result of Ivanov and Turaev \cite{IvTu82}
that the sup norm of the Euler class $\epsilon$ of flat
$\mathrm{SL}(4,\mathbb{R})$-bundles satisfies the inequality
$\left\Vert \epsilon \right\Vert_\infty \leq 1/2^4$, the following
consequence of Theorem \ref{Thm: simpl vol mfld cov by HxH} is
immediate:
\begin{corollary}\label{Cor: Milnor-Wood}Let $M$ be a closed, oriented Riemannian
manifold whose universal cover $\widetilde{M}$ is isometric to $\mathbb{H}%
^{2}\times\mathbb{H}^{2}$. Let $\xi$ be an
$\mathrm{SL}(4,\mathbb{R})$-bundle over $M$. If $\xi$ admits a flat
structure, then
\[
\chi(\xi)=\langle\epsilon(\xi),[M]\rangle\leq \frac{3}{2^3}\chi(M)\text{.}%
\]
\end{corollary}

Note that if $M$ is a product of hyperbolic surfaces and the bundle
$\xi$ is a product of flat bundles, then the stronger inequality
$\chi(\xi)\leq \frac{1}{2^2}\chi(M) $ follows from Milnor's
celebrated inequality \cite{Mi58}. This is probably the correct
bound in the general case also. However, Corollary \ref{Cor:
Milnor-Wood} is good enough to conclude that the considered
manifolds do not admit an affine structure. Indeed, if this was the
case, the $\mathrm{SL}(4,\mathbb{R})$-bundle associated to the
tangent bundle $TM$ could be endowed with a flat structure, and we
would get the impossible inequality
$\chi(M)=\chi(TM)\leq(3/2^3)\chi(M)$.

Another consequence of Theorem \ref{Thm: simpl vol mfld cov by HxH}
is the following product formula:

\begin{corollary}
\label{Corollary: Product for Surf}Let $M$ and $N$ be closed, oriented
surfaces. Then%
\[
\left\Vert M\times N\right\Vert =\frac{3}{2}\left\Vert M\right\Vert
\cdot\left\Vert N\right\Vert .
\]

\end{corollary}

Since the simplicial volume of a surface $\Sigma_{g}$ of genus
$g\geq1$ is equal to $\left\Vert \Sigma_{g}\right\Vert
=4(g-1)=2|\chi (\Sigma_{g})|$, we obtain, for $g,h\geq1$,
\[
\left\Vert \Sigma_{g}\times\Sigma_{h}\right\Vert =24\cdot(g-1)(h-1),
\]
which gives the first exact value of a nonvanishing simplicial volume for a
manifold not admitting a constant curvature metric.

\begin{proof}
[Proof of Corollary \ref{Corollary: Product for Surf}]If either of
$M$ or $N$ is the $2$-sphere or the $2$-torus, then both sides of
the equality vanish trivially. If $M$ and $N$ are endowed with a
hyperbolic structure, then the proportionality principle for
$2$-dimensional hyperbolic manifolds gives us $\pi\left\Vert
M\right\Vert =\mathrm{Vol}(M)$ and $\pi\left\Vert N\right\Vert
=\mathrm{Vol}(N)$, so that by Theorem
\ref{Thm: simpl vol mfld cov by HxH},%
\[
\left\Vert M\times N\right\Vert =\frac{3}{2\pi^{2}}\mathrm{Vol}(M\times
N)=\frac{3}{2\pi^{2}}\mathrm{Vol}(M)\mathrm{Vol}(N)=\frac{3}{2}\left\Vert
M\right\Vert \cdot\left\Vert N\right\Vert ,
\]
as claimed.
\end{proof}

This is the first instance of an exact product formula for the simplicial
volume. Previously known were the rather elementary inequalities%
\[
\left\Vert M\right\Vert \cdot\left\Vert N\right\Vert \leq\left\Vert M\times
N\right\Vert \leq\binom{m+n}{m}\left\Vert M\right\Vert \cdot\left\Vert
N\right\Vert ,
\]
where $M$ and $N$ are any closed, oriented manifolds of dimension $m$ and $n$
respectively. Furthermore, when $M$ and $N$ are hyperbolic surfaces, the upper
bound of $6\cdot\left\Vert M\right\Vert \cdot\left\Vert N\right\Vert $ was
improved to $3,25\cdot\left\Vert M\right\Vert \cdot\left\Vert N\right\Vert $
by Bowen et al. in \cite{BoEtAl04} by exhibiting explicit triangulations of
products of polygons. The authors also give lower bounds for the minimal
number of simplices in such triangulations, which we improve in \cite{Bu07}
using the explicit cocycle $\Theta$ representing the volume form
$\omega_{\mathbb{H}^{2}\times\mathbb{H}^{2}}$. Since those triangulations
produce fundamental cycles for the fundamental class $[M\times N]$ which in
view of Corollary \ref{Corollary: Product for Surf} have strictly greater
$\ell^{1}$-norm than $\left\Vert M\times N\right\Vert $, this indicates the
existence of triangulations of products of hyperbolic surfaces not arising
from triangulations of fundamental domains of the form of a product of
polygons. (Note however that the lower bound for the simplicial volume given
in \cite{BoEtAl04} is incorrect as it relies on the invalid Lemma 2.7.)

In the same way as a hyperbolic surface $M$ can be covered by precisely
$\left\Vert M\right\Vert $ ideal triangles, one could ask if such a covering
can be found for products of hyperbolic surfaces. In fact, the cocycle
representing $\omega_{\mathbb{H}^{2}\times\mathbb{H}^{2}}$ that we exhibit
gives natural candidates for the building blocks of such an ideal
tessellation, since it takes extremal values on very specific $5$-tuples of
points of $\partial\mathbb{H}^{2}\times\partial\mathbb{H}^{2}$.

The present computations are used in \cite{LoSa07} to give the exact
value of the simplicial volume of Hilbert modular surfaces. Those
are {\it open}, finite volume, $\mathbb{Q}$-rank $1$ manifolds with
universal cover isometric to $\mathbb{H}^{2}\times\mathbb{H}^{2}$.

This paper is structured as follows: In Section \ref{Section: volume form}, we
recall the definition of continuous cohomology, and give an explicit
cocycle\ $\Theta$ representing the volume form $\omega_{\mathbb{H}^{2}%
\times\mathbb{H}^{2}}$ in Proposition \ref{Prop: J(omega)=...}. In Section
\ref{Section: cont bdd coho of H}, we introduce continuous bounded cohomology,
and show how our Main Theorem reduces to computing the norm of our explicit
representative $\Theta$ in the cohomology group $H_{c,b}^{4}(H,\widetilde
{\mathbb{R}})$, where $H$ is the (full) isometry group of $\mathbb{H}%
^{2}\times\mathbb{H}^{2}$ and $\widetilde{\mathbb{R}}$ is the real line
endowed with the action of $H$ given by orientation. In Section
\ref{Section: The norm of Theta}, we compute the norm of $\Theta$ both as a
cocycle in Proposition \ref{Prop: norm of cochain Theta} and as a cohomology
class of $H_{c,b}^{4}(H,\widetilde{\mathbb{R}})$ in Theorem 6. Finally, we
prove that the comparison map $H_{c,b}^{4}(H,\widetilde{\mathbb{R}%
})\rightarrow H_{c}^{4}(H,\widetilde{\mathbb{R}})$ is an isomorphism in
Section \ref{Section: The comparison map}.

\subsubsection*{Acknowledgements}

I am indebted to Nicolas Monod for his useful comments on preliminary versions
of this paper.

\section{The volume form $\omega_{\mathbb{H}^{2}\times\mathbb{H}^{2}}$ in
$H_{c}^{4}(\mathrm{PSL}_{2}\mathbb{R}\times\mathrm{PSL}_{2}\mathbb{R}%
,\mathbb{R})$\label{Section: volume form}}

Let $G$ be a topological group and $E$ a $G$-module. Recall (for example from
\cite{Gui80} or \cite{BoWa00}) that the continuous cohomology $H_{c}^{\ast
}(G,E)$ of $G$ with coefficients in $E$ can be computed as the cohomology of
the cocomplex $C^{\ast}(G,E)^{G}$ endowed with its natural symmetric
coboundary operator, where%
\[
C^{q}(G,E)=\left\{  \left.  f:G^{q+1}\longrightarrow\mathbb{R}\text{
}\right\vert \text{ }f\text{ is alternating, measurable}\right\}  ,
\]
and $C^{q}(G,E)^{G}$ denotes the subspace of $G$-invariant cochains, where the
action of $G$ on $C^{q}(G,E)$ is given by
\[
\left(  g\cdot f\right)  \left(  g_{0},...,g_{q}\right)  =g\cdot f(g^{-1}%
g_{0},...,g^{-1}g_{q}),
\]
for every $(g_{0},...,g_{q})$ in $G^{q+1}$, $f$ in $C^{q}(G,E)$ and $g$ in $G$.

Let now $G$ be a Lie group, $K<G$ a maximal compact subgroup and $X=G/K$ the
associated symmetric space. The Van Est isomorphism
\[
\xymatrix{  \mathcal{J}:A^{\ast}(X,E)^{G}\ar[r]^{\ \ \ \ \cong} &
H_{c}^{\ast}(G,E)}
\]
between the $G$-invariant $E$-valued differential forms on $X$ (where the $G$
-action on $A^{\ast}(G,E)$ is defined analogously to that on $C_{c}^{\ast
}(G,E)$) and the continuous cohomology of $G$ with coefficients in $E$ is both
natural and multiplicative. Note furthermore that Dupont gave it an explicit
description at the cochain level in \cite{Du76}.

\subsection*{The volume form in $\mathbb{H}^{2}$}

Let $\omega_{\mathbb{H}^{2}}$ denote the volume form in $A^{2}(\mathbb{H}%
^{2},\mathbb{R})^{\mathrm{PSL}_{2}\mathbb{R}}$. While $\mathcal{J}%
(\omega_{\mathbb{H}^{2}})\in H_{c}^{2}(\mathrm{PSL}_{2}\mathbb{R},\mathbb{R})$
can, by Dupont's description of the Van Est isomorphism be represented by the
cocycle sending a triple of points $(g_{0},g_{1},g_{2})$ in $\left(
\mathrm{PSL}_{2}\mathbb{R}\right)  ^{3}$ to the signed volume of the geodesic
triangle with vertices $(g_{0}x,g_{1}x,g_{2}x)$, for some fixed point $x$ in
$\mathbb{H}^{2}$, let us now describe another cocycle representing
$\mathcal{J}(\omega_{\mathbb{H}^{2}})$. Define%
\[%
\begin{array}
[c]{rccl}%
\mathrm{Or}: & \left(  S^{1}\right)  ^{3} & \longrightarrow & \mathbb{R}\\
& (\xi_{0},\xi_{1},\xi_{2}) & \longmapsto & \left\{
\begin{array}
[c]{rl}%
+1 & \text{ \ \ if }\xi_{0},\xi_{1},\xi_{2}\text{ are positively oriented,}\\
-1 & \text{ \ \ if }\xi_{0},\xi_{1},\xi_{2}\text{ are negatively oriented,}\\
0 & \text{ \ \ if }\xi_{i}=\xi_{j}\text{ for }i\neq j.
\end{array}
\right.
\end{array}
\]
Fix a point $\xi$ in $S^{1}$ and let $\mathrm{Or}_{\xi}:\left(  \mathrm{PSL}%
_{2}\mathbb{R}\right)  ^{3}\rightarrow\mathbb{R}$ be the cocycle defined by%
\[
\mathrm{Or}_{\xi}(g_{0},g_{1},g_{2})=\mathrm{Or}(g_{0}\xi,g_{1}\xi,g_{2}\xi).
\]
It is well known and easy to check that%
\[
\mathcal{J}(\omega_{\mathbb{H}^{2}})=\pi\lbrack\mathrm{Or}_{\xi}].
\]
Thus, the cocycle $\pi\mathrm{Or}_{\xi}$ representing $\mathcal{J}%
(\omega_{\mathbb{H}^{2}})$ can be thought of as sending a triple of points
$(g_{0},g_{1},g_{2})$ to the signed volume of the \textit{ideal} geodesic
triangle with vertices $(g_{0}\xi,g_{1}\xi,g_{2}\xi)$. It is the limit of the
above described cocycle when $x\in\mathbb{H}^{2}$ tends to $\xi\in
\partial\mathbb{H}^{2}=S^{1}$.

\subsection*{The volume form in $\mathbb{H}^{2}\times\mathbb{H}^{2}$}

Let $\omega_{\mathbb{H}^{2}\times\mathbb{H}^{2}}$ denote the volume form in
$A^{4}(\mathbb{H}^{2}\times\mathbb{H}^{2},\mathbb{R})^{\mathrm{PSL}%
_{2}\mathbb{R\times}\mathrm{PSL}_{2}\mathbb{R}}$. We will abuse notation and
write $p_{i}:Y\times Y\rightarrow Y$, for $i=1,2$, for the projections on the
first and second factors for $Y=\mathbb{H}^{2}$, $Y=\partial\mathbb{H}^{2}$ or
$Y=\mathrm{PSL}_{2}\mathbb{R}$. Which of those spaces is meant should be clear
from the context. For $i=1,2$, set%
\begin{align*}
\omega_{i}  &  =p_{i}^{\ast}(\omega_{\mathbb{H}^{2}})\in A^{2}(\mathbb{H}%
^{2}\times\mathbb{H}^{2},\mathbb{R})^{\mathrm{PSL}_{2}\mathbb{R\times
}\mathrm{PSL}_{2}\mathbb{R}}\text{ \ \ and}\\
\mathrm{Or}_{i}  &  =p_{i}^{\ast}(\mathrm{Or}):(S^{1}\times S^{1}%
)^{3}\longrightarrow\mathbb{R}\text{.}%
\end{align*}
Recall that the standard cup product $f_{1}\cup f_{2}$ of a $p$-cochain
$f_{1}:Y^{p+1}\rightarrow\mathbb{R}$ and a $q$-cochain $f_{2}:Y^{q+1}%
\rightarrow\mathbb{R}$ is the \textit{nonalternating }$(p+q)$-cochain sending
the $(p+q+1)$-tuple $(y_{0},...,y_{p+q})$ to the product $f_{1}(y_{0}%
,...,y_{p})\cdot f_{2}(y_{p},...,y_{p+q})$. Moreover, given a non necessarily
alternating $p$-cochain $f:Y^{p+1}\rightarrow\mathbb{R}$, its alternation is
the alternating cochain $\mathrm{Alt}(f):Y^{p+1}\rightarrow\mathbb{R}$ defined
by%
\[
\mathrm{Alt}(f)(y_{0},...,y_{p})=\frac{1}{\left(  p+1\right)  !}\sum
_{\sigma\in\mathrm{Sym}(p+1)}\mathrm{sign}(\sigma)f(y_{\sigma(0)}%
,...,y_{\sigma(p)}),
\]
for every $(y_{0},...,y_{p})$ in $Y^{p+1}$. Set%
\[
\Theta=\mathrm{Alt}(\mathrm{Or}_{1}\cup\mathrm{Or}_{2}):(S^{1}\times
S^{1})^{5}\longrightarrow\mathbb{R}.
\]

Fix a point $\xi$ in $S^{1}$ and let $\Theta_{\xi}:(\mathrm{PSL}%
_{2}\mathbb{R\times}\mathrm{PSL}_{2}\mathbb{R})^{5}\rightarrow\mathbb{R}$ be
the cocycle defined by
\[
\Theta_{\xi}((g_{0},h_{0}),...,(g_{4},h_{4}))=\Theta((g_{0}\xi,h_{0}%
\xi),...,(g_{4}\xi,h_{4}\xi)),
\]
for every $((g_{0},h_{0}),...,(g_{4},h_{4}))$ in $(\mathrm{PSL}_{2}%
\mathbb{R\times}\mathrm{PSL}_{2}\mathbb{R})^{5}$. Note that by construction,%
\[
\Theta_{\xi}=\mathrm{Alt}(p_{1}^{\ast}(\mathrm{Or}_{\xi})\cup p_{2}^{\ast
}(\mathrm{Or}_{\xi})).
\]

\begin{proposition}
\label{Prop: J(omega)=...}$\mathcal{J}(\omega_{\mathbb{H}^{2}\times
\mathbb{H}^{2}})=\pi^{2}\left[  \Theta_{\xi}\right]  .$
\end{proposition}

\begin{proof}
By definition of the Riemannian product metric, the volume form on the
product, is the wedge product of the volume forms on the factors, so that
$\omega_{\mathbb{H}^{2}\times\mathbb{H}^{2}}=\omega_{1}\wedge\omega_{2}$.
Since the van Est isomorphism $\mathcal{J}$\ is multiplicative, we thus have
\[
\mathcal{J}(\omega_{\mathbb{H}^{2}\times\mathbb{H}^{2}})=\mathcal{J}%
(\omega_{1}\wedge\omega_{2})=\mathcal{J}(\omega_{1})\cup\mathcal{J}(\omega
_{2}).
\]
But by naturality of $\mathcal{J}$, we now obtain, for $i=1,2,$%
\[
\mathcal{J}(\omega_{i})=\mathcal{J}(p_{i}^{\ast}(\omega_{\mathbb{H}^{2}%
}))=p_{i}^{\ast}(\mathcal{J}(\omega_{\mathbb{H}^{2}}))=p_{i}^{\ast}(\pi
\lbrack\mathrm{Or}_{\xi}]).
\]
In particular, we get%
\[
\mathcal{J}(\omega_{\mathbb{H}^{2}\times\mathbb{H}^{2}})=\pi^{2}[p_{1}^{\ast
}(\mathrm{Or}_{\xi})]\cup\lbrack p_{2}^{\ast}(\mathrm{Or}_{\xi})]=\pi
^{2}\left[  \Theta_{\xi}\right]  ,
\]
since the cup product is given, at the cochain level, by alternating the
standard cup product.
\end{proof}

\section{The continuous (bounded) cohomology of $H=\mathrm{Isom}%
(\mathbb{H}^{2}\times\mathbb{H}^{2})$ \label{Section: cont bdd coho of H}}

For more details on continuous, bounded cohomology, we invite the reader to
consult \cite{Mo01}. Let $G$ be a topological group and $E$ a Banach
$G$-module. Recall that the continuous cohomology of $G$ with coefficients in
$E$ was defined in the previous section as the cohomology of the cocomplex
$C_{c}^{\ast}(G,E)^{G}$. Now that $E$ is moreover assumed to be a Banach
space, with norm $\left\Vert -\right\Vert _{E}$, say, we can consider the sup
norm%
\[
\left\Vert f\right\Vert _{\infty}=\sup\left\{  \left.  \left\Vert
f(g_{0},...,g_{q})\right\Vert _{E}\text{ }\right\vert \text{ }(g_{0}%
,...,g_{q})\in G^{q+1}\right\}
\]
of any cochain $f$ in $C^{q}(G,E)^{G}$. Clearly, the coboundary operator
restricts to the cocomplex $C_{b}^{\ast}(G,E)^{G}$ of bounded $G$-invariant
cochains, where
\[
C_{b}^{q}(G,E)=\left\{  f\in C_{c}^{q}(G)\left\vert \text{ }\left\Vert
f\right\Vert _{\infty}<+\infty\right.  \right\}  ,
\]
and the \textit{continuous bounded cohomology }$H_{c,b}^{\ast}(G,E)$
\textit{of }$G$\textit{ with coefficients in }$E$ is defined as the cohomology
of this cocomplex. The inclusion of cocomplexes $C_{c,b}^{\ast}(G,E)^{G}%
\subset C_{c}^{\ast}(G,E)^{G}$ induces a \textit{comparison map}
$c:H_{c,b}^{\ast}(G,E)\rightarrow H_{c}^{\ast}(G,E)$. The sup norm defines
both a seminorm on $H_{c}^{\ast}(G,E)$ and $H_{c,b}^{\ast}(G,E)$ and we
continue to denote those by $\left\Vert -\right\Vert _{\infty}$. (Note that on
$H_{c}^{\ast}(G,E)$ we allow the value $+\infty$.) We will abuse terminology
and refer to those seminorms as (sup) norms. By definition, we have for any
$\alpha$ in $H_{c}^{q}(G,E)$:%
\[
\left\Vert \alpha\right\Vert _{\infty}=\inf\left\{  \left.  \left\Vert
\alpha_{b}\right\Vert _{\infty}\text{ }\right\vert \text{ }\alpha_{b}\in
H_{c,b}^{q}(G,E),\text{ }c(\alpha_{b})=\alpha\right\}  ,
\]
where the right hand side of the above equation is understood to be equal to
infinity when the infimum is taken over the empty set.

Let now $G$ be a Lie group, $K<G$ a maximal compact subgroup and $X=G/K$ the
associated symmetric space. We have already come across, in the previous
section, a very convenient cocomplex for the computation of $H_{c}^{\ast
}(G,E)$, namely the degenerate cocomplex $A^{\ast}(X,E)^{G}$ (its differential
is zero since $G$-invariant forms are always closed). Let us now describe
another useful cocomplex for both continuous and continuous, bounded
cohomology: Define%
\[
C_{c}^{q}(X,E)=\left\{  \left.  f:X^{q+1}\rightarrow E\text{ }\right\vert
\text{ }f\text{ continuous, alternating}\right\}
\]
and as above, let $C_{c,b}^{q}(X,E)$ denote its subspace of bounded cochains,
that is the subspace of $C_{c}^{q}(X,E)$ consisting of elements with finite
sup norm. The coboundary operator on $C_{c}^{\ast}(X,E)$, which is the
canonical symmetric operator, clearly restricts to $C_{c,b}^{\ast}(X,E)$. The
action of $G$ on $C_{c}^{\ast}(X,E)$ is defined analogously to the one on
$C^{\ast}(G,E)$. It is a standard fact, that the continuous cohomology
$H_{c}^{\ast}(G,E)$ of $G$\textit{ }with coefficients in $E$ is isomorphic to
the cohomology of the cocomplex $\left(  C_{c}^{\ast}(X,E)^{G},\delta\right)
$ (see for example \cite[Chapitre III, Proposition 2.3]{Gui80} for a proof).
As for the bounded case, it is proven in \cite[Corollary 7.4.10]{Mo01} that
the continuous bounded cohomology $H_{c,b}^{\ast}(G,E)$ of $G$\textit{ }with
coefficients in $E$ is isomorphic to the cohomology of the cocomplex $\left(
C_{c,b}^{\ast}(X,E)^{G},\delta\right)  $. Furthermore, the comparison map
$c:H_{c,b}^{\ast}(G,E)\rightarrow H_{c}^{\ast}(G,E)$ is induced by the natural
inclusion $C_{c,b}^{\ast}(X,E)\subset C_{c}^{\ast}(X,E)$ and the sup norm of
$C_{c}^{\ast}(X,E)$ induces the same seminorms on the cohomology groups
$H_{c}^{\ast}(G,E)$ and $H_{c,b}^{\ast}(G,E)$.

Set $G=\mathrm{PSL}_{2}\mathbb{R}\times\mathrm{PSL}_{2}\mathbb{R}$ and
$H=\mathrm{Isom}(\mathbb{H}^{2}\times\mathbb{H}^{2})$ and note that $G$ is a
subgroup of $H$ of index $8$. Indeed, $\mathrm{PSL}_{2}\mathbb{R}$ has index
$2$ in $\mathrm{Isom}(\mathbb{H}^{2})$, so that $G$ has index $4$ in the
product $\mathrm{Isom}(\mathbb{H}^{2})\times\mathrm{Isom}(\mathbb{H}^{2})$ and
the latter group together with the isometry $\tau$ of $\mathbb{H}^{2}%
\times\mathbb{H}^{2}$ permuting the factors (i.e. $\tau(x,y)=(y,x)$ for
$(x,y)$ in $\mathbb{H}^{2}\times\mathbb{H}^{2}$) generate $H$. Consider the
absolute value norm on $\mathbb{R}$ and denote by $\widetilde{\mathbb{R}}$ the
Banach space $\mathbb{R}$ endowed with the following action of $H$: an element
$h$ of $H$ acts by multiplication by $+1$, respectively $-1$, if $h$
preserves, resp. reverses, the orientation in $\mathbb{H}^{2}\times
\mathbb{H}^{2}$. Observe that restricted to $G$, this action is trivial, and
we denote by $\mathbb{R}$, the Banach space $\mathbb{R}$ endowed with the
trivial action of $G$.

The subgroup inclusion $G<H$ induces maps $i:H_{c}^{\ast}(H,\widetilde
{\mathbb{R}})\rightarrow H_{c}^{\ast}(G,\mathbb{R})$ and $i_{b}:H_{c,b}^{\ast
}(H,\widetilde{\mathbb{R}})\rightarrow H_{c,b}^{\ast}(G,\mathbb{R})$ which are
realized, at the cochain level, by the canonical inclusions of cocomplexes
$C_{c}^{\ast}(\mathbb{H}^{2}\times\mathbb{H}^{2},\widetilde{\mathbb{R}}%
)^{H}\subset C_{c}^{\ast}(\mathbb{H}^{2}\times\mathbb{H}^{2},\mathbb{R})^{G}$
and $C_{c,b}^{\ast}(\mathbb{H}^{2}\times\mathbb{H}^{2},\widetilde{\mathbb{R}%
})^{H}\subset C_{c,b}^{\ast}(\mathbb{H}^{2}\times\mathbb{H}^{2},\mathbb{R}%
)^{G}$ respectively. In particular, both $i$ and $i_{b}$ can not increase norms.

Averaging the value of a $G$-invariant cocycle on a fundamental
domain for $H/G$, it is readily seen that both $i$ and $i_{b}$ admit
left inverses. In fact, those transfer maps can be described
explicitly as follows: Fix an orientation reversing isometry
$\sigma$ of $\mathbb{H}^{2}$ and let $\sigma_{1}$ (respectively
$\sigma_{2}$) be the orientation reversing isometry of
$\mathbb{H}^{2}\times\mathbb{H}^{2}$ acting as $\sigma$ (resp. the
identity) on the first factor and the identity (resp. $\sigma$) on
the second factor. As above, let $\tau\in H$ be the orientation
preserving isometry
permuting the two factors in $\mathbb{H}^{2}\times\mathbb{H}^{2}$. Define%
\[
m:C_{c}^{q}(\mathbb{H}^{2}\times\mathbb{H}^{2},\mathbb{R})^{G}\longrightarrow
C_{c}^{\ast}(\mathbb{H}^{2}\times\mathbb{H}^{2},\widetilde{\mathbb{R}})^{H}%
\]
as%
\begin{align}
m(f)(z_{0},...,z_{q}) &  =\frac{1}{8}\left[  f(z_{0},...,z_{q})-f(\sigma
_{1}z_{0},...,\sigma_{1}z_{q})\right.  \label{Equ: def of m(f)}\\
&  \text{ \ \ \ \ \ \ \ \ }-f(\sigma_{2}z_{0},...,\sigma_{2}z_{q}%
)+f(\sigma_{1}\sigma_{2}z_{0},...,\sigma_{1}\sigma_{2}z_{q})\nonumber\\
&  \text{ \ \ \ \ \ \ \ \ }+f(\tau z_{0},...,\tau z_{q})-f(\sigma_{1}\tau
z_{0},...,\sigma_{1}\tau z_{q})\nonumber\\
&  \left.  \quad\quad\ \ -f(\sigma_{2}\tau z_{0},...,\sigma_{2}\tau
z_{q})+f(\sigma_{1}\sigma_{2}\tau z_{0},...,\sigma_{1}\sigma_{2}\tau
z_{q})\right]  ,\nonumber
\end{align}
for every $f$ in $C^{q}(\mathbb{H}^{2}\times\mathbb{H}^{2},\mathbb{R})^{G}$
and $(z_{0},...,z_{q})$ in $\left(  \mathbb{H}^{2}\times\mathbb{H}^{2}\right)
^{q+1}$. To check that $m$ is well defined, we need to verify that $m(f)$ is
$H$-invariant whenever $f$ is $G$-invariant: First, note that since $\tau$ has
order $2$, $m(f)$ is invariant with respect to $\tau$. Second, we compute%
\begin{align*}
m(f)(\sigma_{1}z_{0},...,\sigma_{1}z_{q}) &  =\frac{1}{8}\left[  f(\sigma
_{1}z_{0},...,\sigma_{1}z_{q})-f(\sigma_{1}^{2}z_{0},...,\sigma_{1}^{2}%
z_{q})\right.  \\
&  \quad\quad\text{\ }-f(\sigma_{2}\sigma_{1}z_{0},...,\sigma_{2}\sigma
_{1}z_{q})+f(\sigma_{1}\sigma_{2}\sigma_{1}z_{0},...,\sigma_{1}\sigma
_{2}\sigma_{1}z_{q})\\
&  \text{ \ \ \ \ \ \ }+f(\tau\sigma_{1}z_{0},...,\tau\sigma_{1}%
z_{q})-f(\sigma_{1}\tau\sigma_{1}z_{0},...,\sigma_{1}\tau\sigma_{1}z_{q})\\
&  \left.  \quad\quad-f(\sigma_{2}\tau\sigma_{1}z_{0},...,\sigma_{2}\tau
\sigma_{1}z_{q})+f(\sigma_{1}\sigma_{2}\tau\sigma_{1}z_{0},...,\sigma
_{1}\sigma_{2}\tau\sigma_{1}z_{q})\right]  .
\end{align*}
Using the facts that $\sigma_{1}$ commutes with $\sigma_{2}$, that $\sigma
_{1}\tau=\tau\sigma_{2}$, that both $\sigma_{1}^{2}$ and $\sigma_{2}^{2}$
belong to $G$ and that $f$ is $G$-invariant, we have%
\begin{align*}
f(\sigma_{1}^{2}z_{0},...,\sigma_{1}^{2}z_{q}) &  =f(z_{0},...,z_{q}),\\
f(\sigma_{2}\sigma_{1}z_{0},...,\sigma_{2}\sigma_{1}z_{q}) &  =f(\sigma
_{1}\sigma_{2}z_{0},...,\sigma_{1}\sigma_{2}z_{q}),\\
f(\sigma_{1}\sigma_{2}\sigma_{1}z_{0},...,\sigma_{1}\sigma_{2}\sigma_{1}z_{q})
&  =f(\sigma_{2}z_{0},...,\sigma_{2}z_{q}),\\
f(\tau\sigma_{1}z_{0},...,\tau\sigma_{1}z_{q}) &  =f(\sigma_{2}\tau
z_{0},...,\sigma_{2}\tau z_{q}),\\
f(\sigma_{1}\tau\sigma_{1}z_{0},...,\sigma_{1}\tau\sigma_{1}z_{q}) &
=f(\sigma_{1}\sigma_{2}\tau z_{0},...,\sigma_{1}\sigma_{2}\tau z_{q}),\\
f(\sigma_{2}\tau\sigma_{1}z_{0},...,\sigma_{2}\tau\sigma_{1}z_{q}) &  =f(\tau
z_{0},...,\tau z_{q}),\\
f(\sigma_{1}\sigma_{2}\tau\sigma_{1}z_{0},...,\sigma_{1}\sigma_{2}\tau
\sigma_{1}z_{q}) &  =f(\sigma_{1}\tau z_{0},...,\sigma_{1}\tau z_{q}).
\end{align*}
Hence, the above expression for $m(f)(\sigma_{1}z_{0},...,\sigma_{1}z_{q})$ is
equal to $-m(f)(z_{0},...,z_{q})$, which proves the invariance of $m(f)$ with
respect to $\sigma_{1}$. The invariance with respect to $\sigma_{2}$ is proven
symmetrically. Third, let $g$ be an isometry in $G$ and observe that since $G$
is normal in $H$, there exists $g_{1},g_{2}$ and $g_{3}$ in $G$ such that
$\sigma_{1}g=g_{1}\sigma_{1}$, $\sigma_{2}g=g_{2}\sigma_{2}$ and $\sigma
_{1}\sigma_{2}g=g_{3}\sigma_{1}\sigma_{2}$. For $k$ in $G$, define
$\overline{k}$ in $G$ as $\overline{k}=\tau k\tau$. Note that if
$k=(k_{1},k_{2})\in\mathrm{PSL}_{2}\mathbb{R\times}\mathrm{PSL}_{2}\mathbb{R}%
$, then $\overline{k}=(k_{2},k_{1})$. We now have%
\begin{align*}
m(f)(gz_{0},...,gz_{q}) &  =\frac{1}{4}\left[  f(gz_{0},...,gz_{q}%
)-f(\sigma_{1}gz_{0},...,\sigma_{1}gz_{q})\right.  \\
&  \text{ \ \ \ \ \ \ \ }-f(\sigma_{2}gz_{0},...,\sigma_{2}gz_{q}%
)+f(\sigma_{1}\sigma_{2}gz_{0},...,\sigma_{1}\sigma_{2}gz_{q})\\
&  \text{ \ \ \ \ \ \ \ \ }+f(\tau gz_{0},...,\tau gz_{q})-f(\sigma_{1}\tau
gz_{0},...,\sigma_{1}\tau gz_{q})\\
&  \left.  \quad\ \ \ \quad-f(\sigma_{2}\tau gz_{0},...,\sigma_{2}\tau
gz_{q})+f(\sigma_{1}\sigma_{2}\tau gz_{0},...,\sigma_{1}\sigma_{2}\tau
gz_{q})\right]  \\
&  =\frac{1}{4}\left[  f(gz_{0},...,gz_{q})-f(g_{1}\sigma_{1}z_{0}%
,...,g_{1}\sigma_{1}z_{q})\right.  \\
&  \text{ \ \ \ \ \ \ \ }-f(g_{2}\sigma_{2}z_{0},...,g_{2}\sigma_{2}%
z_{q})+f(g_{3}\sigma_{1}\sigma_{2}z_{0},...,g_{3}\sigma_{1}\sigma_{2}z_{q})\\
&  \text{ \ \ \ \ \ \ \ \ }+f(\overline{g}\tau z_{0},...,\overline{g}\tau
z_{q})-f(\overline{g}_{2}\sigma_{1}\tau z_{0},...,\overline{g}_{2}\sigma
_{1}\tau z_{q})\\
&  \left.  \quad\ \ \ \quad-f(\overline{g}_{1}\sigma_{2}\tau z_{0}%
,...,\overline{g}_{1}\sigma_{2}\tau z_{q})+f(\overline{g}_{3}\sigma_{1}%
\sigma_{2}\tau z_{0},...,\overline{g}_{3}\sigma_{1}\sigma_{2}\tau
z_{q})\right]  \\
&  =m(f)(z_{0},...,z_{q}),
\end{align*}
where for the last equality we have used eight times the $G$-invariance of
$f$. Finally, the $H$-invariance of $m(f)$ follows from that $H$ is generated
by $\sigma_{1},\sigma_{2},\tau$ and $G$. Observe also that, by the
$G$-invariance of $f$, the definition of $m(f)$ is independent of the choice
of $\sigma$.

It is readily seen that $m$ is a cochain map which moreover restricts to a map%
\[
m_{b}:C_{c,b}^{q}(\mathbb{H}^{2}\times\mathbb{H}^{2},\mathbb{R})^{G}%
\longrightarrow C_{c,b}^{q}(\mathbb{H}^{2}\times\mathbb{H}^{2},\widetilde
{\mathbb{R}})^{H}%
\]
between the respective bounded cocomplexes. In particular, $m$ and $m_{b}$
induce maps, which we still denote by $m$ and $m_{b}$ between the
corresponding cohomology groups. It is clear that neither $m$ nor $m_{b}$ can
increase norms. Furthermore, since both the inclusion of cocomplexes
$C_{c}^{\ast}(\mathbb{H}^{2}\times\mathbb{H}^{2},\widetilde{\mathbb{R}}%
)^{H}\subset C_{c}^{\ast}(\mathbb{H}^{2}\times\mathbb{H}^{2},\mathbb{R})^{G}$
composed with $m$ and $C_{c,b}^{\ast}(\mathbb{H}^{2}\times\mathbb{H}%
^{2},\widetilde{\mathbb{R}})^{H}\subset C_{c,b}^{\ast}(\mathbb{H}^{2}%
\times\mathbb{H}^{2},\mathbb{R})^{G}$ composed with $m_{b}$ are the identity
maps, we have obtained a commutative diagram%

\[
\xymatrix{ H^{\ast}_c(H,\widetilde{\mathbb{R}})
\ar@/^1.5pc/[rr]_{\mathrm{Id}}\ar@{^{(}->}[r]^<<<<<{i} &
H_c^*(G,\mathbb{R})
\ar@{->>}[r]^<<<<<{m} &H^{\ast}_c(H,\widetilde{\mathbb{R}}) \\
H^{\ast}_ {c,b}(H,\widetilde{\mathbb{R}})
\ar@/_1.5pc/[rr]^{\mathrm{Id}}\ar[u]_c \ar@{^{(}->}[r]^<<<<<{i_b} &
H_{c,b}^*(G,\mathbb{R}) \ar[u]_c \ar@{->>}[r]^<<<<<{m_b} &
H^{\ast}_{c,b}(H,\widetilde{\mathbb{R}}). \ar[u]_c }
\]

Going back to the definition of the above cohomology groups in terms of the
cocomplexes $C^{\ast}(G,\mathbb{R})^{G}$ and $C^{\ast}(H,\widetilde
{\mathbb{R}})^{H}$ and their bounded subcocomplexes, we see that, for any
$\xi$ in $S^{1}$, the cochain $\Theta_{\xi}\in C^{4}(G,\mathbb{R})^{G}$
encountered in the previous section has sup norm $\left\Vert \Theta_{\xi
}\right\Vert _{\infty}\leq1$ and hence belongs to $C_{b}^{4}(G,\mathbb{R}%
)^{G}$. (In fact, we will show in Proposition
\ref{Prop: norm of cochain Theta} that $\left\Vert \Theta_{\xi}\right\Vert
_{\infty}=\left\Vert \Theta\right\Vert _{\infty}=2/3$.) Denote by $\left[
\Theta\right]  \in H_{c}^{4}(G,\mathbb{R})$ and $\left[  \Theta\right]
_{b}\in H_{c,b}^{4}(G,\mathbb{R})$ the corresponding cohomology classes.

Since $H$ acts on $S^{1}\times S^{1}=\partial\mathbb{H}^{2}\times
\partial\mathbb{H}^{2}$, we can as well extend our definition of $\Theta_{\xi
}$ to a cocycle on $H$ which maps a $5$-tuple $(h_{0},...,h_{4})$ in $H^{5}$
to $\Theta(h_{0}(\xi,\xi),...,h_{4}(\xi,\xi))$. The extended $\Theta_{\xi}$ is
clearly $H$-invariant and hence belongs to $C^{4}(H,\widetilde{\mathbb{R}%
})^{H}$ and also to $C_{b}^{4}(H,\widetilde{\mathbb{R}})^{H}$. Denote by
$\left[  \Theta\right]  ^{H}\in H_{c}^{4}(H,\widetilde{\mathbb{R}})$ and
$\left[  \Theta\right]  _{b}^{H}\in H_{c,b}^{4}(H,\widetilde{\mathbb{R}})$ the
corresponding cohomology classes.

Because the four cohomology classes $\left[  \Theta\right]  ,\left[
\Theta\right]  ^{H},\left[  \Theta\right]  _{b}$ and $\left[  \Theta\right]
_{b}^{H}$ are the equivalence classes of the same cocycle $\Theta_{\xi}$, it
is obvious that
\[%
\begin{array}
[c]{ll}%
c(\left[  \Theta\right]  _{b})=\left[  \Theta\right]  ,\text{ \ \ \ \ \ \ } &
c(\left[  \Theta\right]  _{b}^{H})=\left[  \Theta\right]  ^{H}\text{
\ \ and}\\
i(\left[  \Theta\right]  ^{H})=\left[  \Theta\right]  , & i(\left[
\Theta\right]  _{b}^{H})=\left[  \Theta\right]  _{b}.
\end{array}
\]

\begin{proposition}
\label{Prop: [Theta]=[Theta]H}$\left\Vert \left[  \Theta\right]  \right\Vert
_{\infty}=\left\Vert \left[  \Theta\right]  ^{H}\right\Vert _{\infty}.$
\end{proposition}

\begin{proof}
We use the facts that both $i$ and $m$ can not increase norms and that $m\circ
i$ is the identity on $H_{c}^{4}(H,\widetilde{\mathbb{R}})$ to obtain%
\[
\left\Vert \left[  \Theta\right]  ^{H}\right\Vert _{\infty}=\left\Vert \left(
m\circ i\right)  \left(  \left[  \Theta\right]  ^{H}\right)  \right\Vert
_{\infty}\leq\left\Vert i\left(  \left[  \Theta\right]  ^{H}\right)
\right\Vert _{\infty}\leq\left\Vert \left[  \Theta\right]  ^{H}\right\Vert
_{\infty}.
\]
In particular, all the above inequalities are equalities. The proposition is
now immediate from that $i\left(  \left[  \Theta\right]  ^{H}\right)  =\left[
\Theta\right]  $.
\end{proof}

\begin{theorem}
\label{Thm: comparison map}The comparison map $c:H_{cb}^{\ast}(H,\widetilde
{\mathbb{R}})\rightarrow H_{c}^{\ast}(H,\widetilde{\mathbb{R}})$ is an
isomorphism in degree $4$ sending $\left[  \Theta\right]  _{b}^{H}$ to
$\left[  \Theta\right]  ^{H}$.
\end{theorem}

Note that if the comparison map is an isomorphism, it is automatically
isometric. Theorem \ref{Thm: comparison map} is proven in Section
\ref{Section: The comparison map} in a constructive way: given a cocycle
representing a cohomology class of degree $4$ in the kernel of the comparison
map, we can explicitly express it as a coboundary. This isomorphism allows us
to reduce the proof of our Main Theorem to the computation of the norm of
$\Theta$ in $H_{cb}^{\ast}(H,\widetilde{\mathbb{R}})$, as stated in the next
theorem, which we will prove in Section \ref{Section: The norm of Theta}.
\setcounter{theorem}{6}


\begin{thmForTheta}
\label{Thm: norm coho class Theta H b}The norm of $\Theta$ in $H_{cb}^{\ast
}(H,\widetilde{\mathbb{R}})$ is equal to $\left\Vert \left[  \Theta\right]
_{b}^{H}\right\Vert _{\infty}=2/3$.

\end{thmForTheta}

\begin{proof}
[Proof of Main Theorem]Recall that we showed in Proposition
\ref{Prop: J(omega)=...} that $\omega_{\mathbb{H}^{2}\times\mathbb{H}^{2}}\in
H_{c}^{4}(G,\mathbb{R)}$ (which was there denoted by $\mathcal{J}%
(\omega_{\mathbb{H}^{2}\times\mathbb{H}^{2}})$) is equal to $\pi^{2}%
[\Theta_{\xi}]=\pi^{2}[\Theta]$. Applying successively Proposition
\ref{Prop: [Theta]=[Theta]H}, Theorem \ref{Thm: comparison map} and Theorem 6,
we obtain
\[
\frac{\left\Vert \omega_{\mathbb{H}^{2}\times\mathbb{H}^{2}}\right\Vert }%
{\pi^{2}}=\left\Vert \left[  \Theta\right]  \right\Vert _{\infty}=\left\Vert
\left[  \Theta\right]  ^{H}\right\Vert _{\infty}=\left\Vert \left[
\Theta\right]  _{b}^{H}\right\Vert _{\infty}=2/3.
\]

\end{proof}

We now introduce yet another cocomplex for the computation of the continuous,
bounded cohomology groups: Let $C_{b}^{q}(S^{1}\times S^{1},\widetilde
{\mathbb{R}})$ denote the space of alternating, measurable, bounded,
real-valued functions on $(S^{1}\times S^{1})^{q+1}$ endowed with its natural
symmetric coboundary operator $\delta$. The action of $H$ on $C_{b}^{q}%
(S^{1}\times S^{1},\widetilde{\mathbb{R}})$ is defined analogously to the one
of $H$ on $C^{q}(H,\widetilde{\mathbb{R}})$. Its subspaces of $H$-invariant
and $G$-invariant functions are denoted by $C_{b}^{q}(S^{1}\times
S^{1},\widetilde{\mathbb{R}})^{H}$ and $C_{b}^{q}(S^{1}\times S^{1}%
,\mathbb{R})^{G}$ respectively. It is proven in \cite[Corollary 7.5.9]{Mo01}
that the cohomology of the latter cocomplexes are isomorphic to $H_{cb}^{\ast
}(H,\widetilde{\mathbb{R}})$ and $H_{cb}^{\ast}(G,\mathbb{R})$ respectively,
that the map $i_{b}:H_{cb}^{\ast}(H,\widetilde{\mathbb{R}})\rightarrow
H_{cb}^{\ast}(G,\mathbb{R})$ is realized at the cochain level by the canonical
inclusion $C_{b}^{q}(S^{1}\times S^{1},\widetilde{\mathbb{R}})^{H}\subset
C_{b}^{q}(S^{1}\times S^{1},\mathbb{R})^{G}$, and that the sup norm on
$C_{b}^{q}(S^{1}\times S^{1},\widetilde{\mathbb{R}})$ gives rise to the
desired seminorms on the continuous, bounded cohomology groups. Furthermore,
it is easy to verify, that the map $m_{b}:H_{cb}^{\ast}(G,\mathbb{R}%
)\rightarrow H_{cb}^{\ast}(H,\widetilde{\mathbb{R}})$ is realized at the
cochain level by the map
\[
m_{b}:C_{b}^{q}(S^{1}\times S^{1},\mathbb{R)}^{G}\longrightarrow C_{b}%
^{q}(S^{1}\times S^{1},\widetilde{\mathbb{R}}\mathbb{)}^{H}%
\]
defined exactly as in $\left(  \ref{Equ: (Or,z)= 4 times orientation cocycle}%
\right)  $ except that $f$ is now taken in $C_{b}^{q}(S^{1}\times
S^{1},\mathbb{R)}^{G}$ and $(z_{0},...,z_{q})$ in $\left(  S^{1}\times
S^{1}\right)  ^{q+1}$. Finally, note that the cohomology classes $\left[
\Theta\right]  _{b}^{H}$ and $\left[  \Theta\right]  _{b}$ are represented in
the cocomplexes $C_{b}^{q}(S^{1}\times S^{1},\widetilde{\mathbb{R}})^{H}$\ and
$C_{b}^{q}(S^{1}\times S^{1},\mathbb{R})^{G}$ by the cocycle $\Theta
:(S^{1}\times S^{1})^{5}\rightarrow\mathbb{R}$.

\section{The norm of $\Theta$\label{Section: The norm of Theta}}

Recall that $\Theta\in C_{b}^{4}(S^{1}\times S^{1},\widetilde{\mathbb{R}}%
)^{H}$ is defined as the cocycle
\[
\Theta=\mathrm{Alt}(\mathrm{Or}_{1}\cup\mathrm{Or}_{2}),
\]
where $\mathrm{Or}_{1}$ and $\mathrm{Or}_{2}$ are the pullbacks under the
first and second projection respectively of the orientation cocycle on $S^{1}$.

\begin{proposition}
\label{Prop: norm of cochain Theta}$\left\Vert \Theta\right\Vert _{\infty
}=2/3$.
\end{proposition}

\begin{proof}
By definition, we have, for any $5$-tuple $((x_{0},y_{0}),...,(x_{4},y_{4}))$
in $(S^{1}\times S^{1})^{5}$, that $\Theta((x_{0},y_{0}),...,(x_{4},y_{4}))$
is equal to
\[
\frac{1}{120}\sum_{\sigma\in\mathrm{Sym}(5)}\mathrm{sign}(\sigma
)\mathrm{Or}(x_{\sigma(0)},x_{\sigma(1)},x_{\sigma(2)})\cdot\mathrm{Or}%
(y_{\sigma(2)},y_{\sigma(3)},y_{\sigma(4)}).
\]
Set $\tau=(0$ $1$ $...$ $4)$ and observe that every permutation $\sigma
\in\mathrm{Sym}(5)$ can be written uniquely as $\sigma=\tau^{k}\circ\alpha$,
where $\alpha\in\mathrm{Sym}(5)$ maps $2$ to $0$, and $k$ is an integer
between $0$ and $4$. Now, exploiting the fact that $\mathrm{Or}$ is
alternating, we can rewrite the above expression for $\Theta((x_{0}%
,y_{0}),...,(x_{4},y_{4}))$ as
\begin{align}
&  \frac{1}{30}\sum_{\substack{\tau=(0\text{ }1...4)^{k}\\k\in\{0,1,...,4\}}%
}\left[  \mathrm{Or}\left(  x_{\tau(0)},x_{\tau(1)},x_{\tau(2)}\right)
\mathrm{\cdot Or}(y_{\tau(0)},y_{\tau(3)},y_{\tau(4)})\right.
\label{Equation: Basic formula for Theta}\\
&  \text{ \ \ \ \ \ \ \ \ \ \ \ \ \ \ \ \ \ \ }+\mathrm{Or}\left(  x_{\tau
(0)},x_{\tau(3)},x_{\tau(4)}\right)  \mathrm{\cdot Or}(y_{\tau(0)},y_{\tau
(1)},y_{\tau(2)})\nonumber\\
&  \text{ \ \ \ \ \ \ \ \ \ \ \ \ \ \ \ \ \ \ }-\mathrm{Or}\left(  x_{\tau
(0)},x_{\tau(1)},x_{\tau(3)}\right)  \mathrm{\cdot Or}(y_{\tau(0)},y_{\tau
(2)},y_{\tau(4)})\nonumber\\
&  \text{ \ \ \ \ \ \ \ \ \ \ \ \ \ \ \ \ \ \ }-\mathrm{Or}\left(  x_{\tau
(0)},x_{\tau(2)},x_{\tau(4)}\right)  \mathrm{\cdot Or}(y_{\tau(0)},y_{\tau
(1)},y_{\tau(3)})\nonumber\\
&  \text{ \ \ \ \ \ \ \ \ \ \ \ \ \ \ \ \ \ \ }+\mathrm{Or}\left(  x_{\tau
(0)},x_{\tau(1)},x_{\tau(4)}\right)  \mathrm{\cdot Or}(y_{\tau(0)},y_{\tau
(2)},y_{\tau(3)})\nonumber\\
&  \left.  \text{ \ \ \ \ \ \ \ \ \ \ \ \ \ \ \ \ \ }+\mathrm{Or}\left(
x_{\tau(0)},x_{\tau(2)},x_{\tau(3)}\right)  \mathrm{\cdot Or}(y_{\tau
(0)},y_{\tau(1)},y_{\tau(4)})\right]  .\nonumber
\end{align}

Let us now compute the absolute value of the evaluation of $\Theta$ on an
arbitrary $5$-tuple $((x_{0},y_{0}),...,(x_{4},y_{4}))$ in $(S^{1}\times
S^{1})^{5}$. If the $x_{i}$'s are all distinct, we can, since $\Theta$ is
alternating, up to permuting the $x_{i}$'s assume that they are cyclically
ordered according to their numbering. Thus, $\mathrm{Or}\left(  x_{i}%
,x_{j},x_{k}\right)  =+1$ whenever $0\leq i<j<k\leq4$, and all the orientation
cocycles involving the $x_{i}$'s in the expression $\left(
\ref{Equation: Basic formula for Theta}\right)  $ are equal to $+1$, so that
$\Theta((x_{0},y_{0}),...,(x_{4},y_{4}))$ is equal to%
\begin{align*}
&  \frac{1}{30}\sum_{\substack{\tau=(0\text{ }1...4)^{k}\\k\in\{0,1,...,4\}}%
}\left[  \mathrm{Or}(y_{\tau(0)},y_{\tau(3)},y_{\tau(4)})+\mathrm{Or}%
(y_{\tau(0)},y_{\tau(1)},y_{\tau(2)})\right. \\
&  \text{ \ \ \ \ \ \ \ \ \ \ \ \ \ \ \ \ \ \ \ }-\mathrm{Or}(y_{\tau
(0)},y_{\tau(2)},y_{\tau(4)})-\mathrm{Or}(y_{\tau(0)},y_{\tau(1)},y_{\tau
(3)})\\
&  \left.  \text{ \ \ \ \ \ \ \ \ \ \ \ \ \ \ \ \ \ \ }+\mathrm{Or}%
(y_{\tau(0)},y_{\tau(2)},y_{\tau(3)})+\mathrm{Or}(y_{\tau(0)},y_{\tau
(1)},y_{\tau(4)})\right] \\
&  =\frac{1}{30}\sum_{\substack{\tau=(0\text{ }1...4)^{k}\\k\in\{0,1,...,4\}}%
}\left[  \mathrm{Or}(y_{\tau(2)},y_{\tau(3)},y_{\tau(4)})+\mathrm{Or}%
(y_{\tau(0)},y_{\tau(1)},y_{\tau(2)})\right. \\
&  \left.  \text{ \ \ \ \ \ \ \ \ \ \ \ \ \ \ \ \ \ \ \ \ \ \ }-\mathrm{Or}%
(y_{\tau(0)},y_{\tau(1)},y_{\tau(3)})+\mathrm{Or}(y_{\tau(0)},y_{\tau
(1)},y_{\tau(4)})\right]  ,
\end{align*}
where we have used the cocycle relation%
\begin{align*}
0  &  =\delta\mathrm{Or}(y_{\tau(0)},y_{\tau(2)},y_{\tau(3)},y_{\tau(4)})\\
&  =\mathrm{Or}(y_{\tau(2)},y_{\tau(3)},y_{\tau(4)})-\mathrm{Or}(y_{\tau
(0)},y_{\tau(3)},y_{\tau(4)})\\
&  \text{ \ \ \ \ }+\mathrm{Or}(y_{\tau(0)},y_{\tau(2)},y_{\tau(4)}%
)-\mathrm{Or}(y_{\tau(0)},y_{\tau(2)},y_{\tau(3)}).
\end{align*}
It is now immediate that $\left\vert \Theta((x_{0},y_{0}),...,(x_{4}%
,y_{4}))\right\vert \leq2/3$, since the last expression for $\Theta
((x_{0},y_{0}),...,(x_{4},y_{4}))$ is a sum of $5\cdot4=20$ elements admitting
the values $\pm1/30$ and $0$.

If the $x_{i}$'s are not all distinct, then we can without loss of generality
assume that $x_{0}=x_{1}$. If the $y_{i}$'s were all distinct, then by
symmetry we could apply the above argument to show that $\left\vert
\Theta((x_{0},y_{0}),...,(x_{4},y_{4}))\right\vert \leq2/3$. Let us thus
assume that the $y_{i}$'s are not all distinct. If $y_{0}=y_{1}$ then
$(x_{0},y_{0})=(x_{1},y_{1})$ and hence $\Theta((x_{0},y_{0}),...,(x_{4}%
,y_{4}))=0$ since $\Theta$ is alternating. We can now, again without loss of
generality assume that $y_{2}=y_{k}$, for $k$ in $\{0,3,4\}$. In the
expression $\left(  \ref{Equation: Basic formula for Theta}\right)  $ there
are exactly $9$ summands which have as a factor $\mathrm{Or}(x_{0},x_{1}%
,x_{j})$, up to permutation of the entries, for $j\geq2$, and hence vanish.
Furthermore, the summand $\mathrm{Or}(x_{0},x_{3},x_{4})\cdot\mathrm{Or}%
(y_{1},y_{2},y_{k})$ (which exists in $\left(
\ref{Equation: Basic formula for Theta}\right)  $ again up to permutation of
the entries) also vanishes and clearly has not yet been counted among the
summands having a factor of the form $\mathrm{Or}(x_{0},x_{1},x_{j})$, so at
least $10$ of the $30$ summands in $\left(
\ref{Equation: Basic formula for Theta}\right)  $ vanish and
\[
\left\vert \Theta((x_{0},y_{0}),...,(x_{4},y_{4}))\right\vert \leq2/3.
\]
We have thus proven $\left\Vert \Theta\right\Vert _{\infty}\leq\frac{2}{3}$.
To prove equality, observe that if $x_{0},...,x_{4}$ in $S^{1}$ are positively
cyclically ordered according to their numbering and $y_{0},...,y_{4}$ in
$S^{1}$ are $\varepsilon$-cyclically ordered, for $\varepsilon$ in
$\{-1,+1\}$, according to their numbering, then%
\[
\Theta((x_{0},y_{0}),(x_{1},y_{2}),(x_{2},y_{4}),(x_{3},y_{1}),(x_{4}%
,y_{3}))=\varepsilon\cdot\frac{2}{3},
\]
which finishes the proof of the proposition.
\end{proof}

Let us now state and prove three easy lemmas which will furthermore be useful
again in the next section. For the moment, they will allow us a better
understanding of the spaces $C_{b}^{\ast}(S^{1}\times S^{1},\widetilde
{\mathbb{R}})^{H}$ in low degree.

\begin{lemma}
\label{Lemma: c(z)=signs c(z)}Let $f$ be a cochain in $C_{b}^{q}(S^{1}\times
S^{1},\widetilde{\mathbb{R}})^{H}$ and let $\underline{z}=((x_{0}%
,y_{0}),...,(x_{q},y_{q}))$ be a $(q+1)$-tuple in $(S^{1}\times S^{1})^{q+1}$.
If there exists $\sigma$ in $\mathrm{Sym}(q+1)$ such that the permutations
$x_{i}\mapsto x_{\sigma(i)}$ and $y_{i}\mapsto y_{\sigma(i)}$, for $0\leq
i\leq q$, can be realized by isometries $g$ and $h$ of $\mathbb{H}^{2}$
respectively, then%
\[
f(\underline{z})=\mathrm{sign}(\sigma)\mathrm{sign}(g)\mathrm{sign}%
(h)f(\underline{z}),
\]
where $\mathrm{sign}(k)=+1$, respectively $-1$, if $k$ is an orientation
preserving, resp. reversing, isometry of $\mathbb{H}^{2}$.
\end{lemma}

\begin{proof}
On the one hand, we have, since $f$ is alternating,
\[
f(\underline{z})=\mathrm{sign}(\sigma)f((x_{\sigma(0)},y_{\sigma
(0)}),...,(x_{\sigma(q)},y_{\sigma(q)})).
\]
On the other hand, using the $H$-invariance of $f$, we get%
\[
f(\underline{z})=\mathrm{sign}(g)\mathrm{sign}(h)f((gx_{0},hy_{0}%
),...,(gx_{q},hy_{q})).
\]
But by assumption, $(x_{\sigma(i)},y_{\sigma(i)})=(gx_{i},hy_{i})$, for every
$0\leq i\leq q$, and the lemma follows.
\end{proof}

\begin{lemma}
\label{Lemma: 4 points, geodesics intersect or not}Let $x_{0},x_{1}%
,x_{2},x_{3}$ be distinct points on $S^{1}$. Denote by $\left\langle
x_{i},x_{j}\right\rangle $, for $i\neq j$, the geodesic in $\mathbb{H}^{2}$
between $x_{i}$ and $x_{j}$ in $\partial\mathbb{H}^{2}=S^{1}$.

\begin{enumerate}
\item If $\left\langle x_{0},x_{1}\right\rangle \cap\left\langle x_{2}%
,x_{3}\right\rangle \neq\varnothing$, then there exists an orientation
preserving isometry of $\mathbb{H}^{2}$ realizing the permutation $(0$ $1)(2$
$3)$.

\item If $\left\langle x_{0},x_{1}\right\rangle \cap\left\langle x_{2}%
,x_{3}\right\rangle =\varnothing$, then there exists an orientation reversing
isometry of $\mathbb{H}^{2}$ realizing the permutation $(0$ $1)(2$ $3)$.
\end{enumerate}
\end{lemma}

\begin{proof}
This is elementary from hyperbolic geometry:

\begin{enumerate}
\item Since the points are all distinct, if the geodesics intersect, they
intersect in precisely one point. Then the rotation by $\pi$ centered at the
intersection realizes the permutation $x_{0}\leftrightarrow x_{1}%
,x_{2}\leftrightarrow x_{3}$ and clearly preserves orientation.

\item There exists a unique geodesic $\gamma$ perpendicular to both
$\left\langle x_{0},x_{1}\right\rangle \ $and $\left\langle x_{2}%
,x_{3}\right\rangle $. The reflection along $\gamma$ is a reversing
orientation isometry of $\mathbb{H}^{2}$ realizing the permutations
$x_{0}\leftrightarrow x_{1},x_{2}\leftrightarrow x_{3}$.
\end{enumerate}
\end{proof}

It is easy to conclude, from Lemma \ref{Lemma: c(z)=signs c(z)}, that
$C_{b}^{q}(S^{1}\times S^{1},\widetilde{\mathbb{R}})^{H}=0$ for $0\leq q\leq
2$, and consequently also $H_{c,b}^{q}(H,\widetilde{\mathbb{R}})=0$ for $0\leq
q\leq2$. In degree $3$, the space of cochain $C_{b}^{3}(S^{1}\times
S^{1},\widetilde{\mathbb{R}})^{H}$ is not zero, but we have the following
useful vanishing criterion:

\begin{lemma}
\label{Lemma: 3-cochain(twisted)=0}Let $f$ be a cochain in $C_{b}^{3}%
(S^{1}\times S^{1},\widetilde{\mathbb{R}})^{H}$ and let $\underline{z}%
=((x_{0},y_{0}),...,(x_{3},y_{3}))$ be a $4$-tuple in $(S^{1}\times S^{1}%
)^{4}$ such that $\left\langle x_{0},x_{2}\right\rangle \cap\left\langle
x_{1},x_{3}\right\rangle \neq\varnothing$ and $\left\langle y_{0}%
,y_{2}\right\rangle \cap\left\langle y_{1},y_{3}\right\rangle =\varnothing$.
Then%
\[
f(\underline{z})=0\text{.}%
\]

\end{lemma}

\begin{proof}
By Lemma \ref{Lemma: 4 points, geodesics intersect or not} there exists an
orientation preserving isometry $g$ of $\mathbb{H}^{2}$ realizing the
permutation $(0$ $2)(1$ $3)$ of the points $x_{0},x_{1},x_{2},x_{3}$ and an
orientation reversing isometry $h$ of $\mathbb{H}^{2}$ realizing the
permutation $(0$ $2)(1$ $3)$ of the points $y_{0},y_{1},y_{2},y_{3}$. Thus, by
Lemma \ref{Lemma: c(z)=signs c(z)} we obtain $f(\underline{z})=-f(\underline
{z})$, and the lemma is proven.
\end{proof}

\begin{thmForTheta}
$\left\Vert [\Theta]_{b}^{H}\right\Vert _{\infty}=2/3$.
\end{thmForTheta}

\begin{proof}
From Proposition \ref{Prop: norm of cochain Theta}, one equality is already
immediate, namely $\left\Vert [\Theta]_{b}^{H}\right\Vert _{\infty}%
\leq\left\Vert \Theta\right\Vert _{\infty}=2/3$. For the other inequality, let
$b\in C_{b}^{3}(S^{1}\times S^{1},\widetilde{\mathbb{R}})^{H}$ be an arbitrary
cochain. As in the end of the proof of Proposition
\ref{Prop: norm of cochain Theta}, let $x_{0},...,x_{4}$, and respectively
$y_{0},...,y_{4}$, be positively cyclically ordered points in $S^{1}$ and
consider the $5$-tuple
\[
((x_{0},y_{0}),(x_{1},y_{2}),(x_{2},y_{4}),(x_{3},y_{1}),(x_{4},y_{3}))
\]
(which we already know has value $2/3$ on $\Theta$). Whatever coordinate one
removes from this given $5$-tuple, the remaining $4$-tuple satisfies the
conditions of Lemma \ref{Lemma: 3-cochain(twisted)=0}, so that%
\[
\delta b((x_{0},y_{0}),(x_{1},y_{2}),(x_{2},y_{4}),(x_{3},y_{1}),(x_{4}%
,y_{3}))=0.
\]
In particular, we obtain
\begin{align*}
\left\Vert \Theta+\delta b\right\Vert _{\infty}  &  \geq\left\vert \left(
\Theta+\delta b\right)  ((x_{0},y_{0}),(x_{1},y_{2}),(x_{2},y_{4}%
),(x_{3},y_{1}),(x_{4},y_{3}))\right\vert \\
&  =\left\vert \Theta((x_{0},y_{0}),(x_{1},y_{2}),(x_{2},y_{4}),(x_{3}%
,y_{1}),(x_{4},y_{3}))\right\vert =2/3,
\end{align*}
and hence%
\[
\left\Vert \lbrack\Theta]_{b}^{H}\right\Vert _{\infty}=\inf\left\{  \left.
\left\Vert \Theta+\delta b\right\Vert _{\infty}\right\vert \text{ }b\in
C_{b}^{3}(S^{1}\times S^{1},\widetilde{\mathbb{R}})^{H}\right\}  \geq2/3,
\]
which finishes the proof of the theorem.
\end{proof}

\section{The comparison map $c:H_{cb}^{\ast}(H,\widetilde{\mathbb{R}%
})\rightarrow H_{c}^{\ast}(H,\widetilde{\mathbb{R}})$%
\label{Section: The comparison map}}

In this last section, we prove Theorem \ref{Thm: comparison map}, that is, we
prove that the comparison map $c:H_{cb}^{\ast}(H,\widetilde{\mathbb{R}%
})\rightarrow H_{c}^{\ast}(H,\widetilde{\mathbb{R}})$ is an isomorphism in
degree $4$ sending $[\Theta]_{b}^{H}$ to $[\Theta]^{H}$. Since $H_{c}^{\ast
}(H,\widetilde{\mathbb{R}})$ injects in the $1$-dimensional cohomology group
$H_{c}^{4}(G,\mathbb{R})$ and contains the nonzero class $[\Theta]^{H}$, it is
clearly also $1$-dimensional, generated by $[\Theta]^{H}$. Moreover, we have
already seen that $c([\Theta]_{b}^{H})=[\Theta]^{H}$. Thus, it only remains to
prove that the comparison map is injective in degree $4$.

Let $p_{i}:S^{1}\times S^{1}\rightarrow S^{1}$, for $i=1,2$, denote the
projection on the first and second factor respectively. For any $(q+1)$-tuple
$(z_{0},...,z_{q})$ in $(S^{1}\times S^{1})^{q+1}$, define%
\[
n_{i}(z_{0},...,z_{q})=\sharp\{p_{i}(z_{0}),...,p_{i}(z_{q})\},
\]
for $i=1,2$. We will now prove inductively on $(n_{1},n_{2})$ that if a
cocycle $f\in C_{b}^{4}(S^{1}\times S^{1},\widetilde{\mathbb{R}})^{H}$
represents a cohomology class which is mapped to zero by the comparison map,
then $f=\delta h$ on $5$-tuples $\underline{z}$ satisfying $n_{1}%
(\underline{z})\leq n_{1}$, $n_{2}(\underline{z})\leq n_{2}$. Observe that the
fact that $c([f])=0$ will only be used in Step \ref{Step 1}, where we show
that $f=\delta h$ on $5$-tuples $\underline{z}$ verifying $n_{1}(\underline
{z})=n_{2}(\underline{z})=3$. Thus, Step \ref{Step 2} and Step \ref{Step 3}
amount to proving that a cocycle vanishing on $5$-tuples $\underline{z}$
satisfying $n_{1}(\underline{z})=n_{2}(\underline{z})=3$ is a coboundary.

\begin{step}
\label{Step 0}Let $f$ be a cochain in $C_{b}^{q}(S^{1}\times S^{1}%
,\widetilde{\mathbb{R}})^{H}$ and $\underline{z}$ be a $(q+1)$-tuple in
$(S^{1}\times S^{1})^{q+1}$. If $n_{1}(\underline{z})\leq2$ or $n_{2}%
(\underline{z})\leq2$, then $f(\underline{z})=0$.
\end{step}

\begin{proof}
By symmetry, it is enough to treat the case $n_{1}(\underline{z})\leq2$. If
$n_{1}(\underline{z})\leq2$, for $\underline{z}=((x_{0},y_{0}),...,(x_{q}%
,y_{q}))$ in $(S^{1}\times S^{1})^{q+1}$, then there exists an orientation
reversing isometry of $\mathbb{H}^{2}$ fixing $x_{0},...,x_{q}$, while the
identity fixes $y_{0},...,y_{q}$. In particular, by Lemma
\ref{Lemma: c(z)=signs c(z)}, $f(\underline{z})=0$.
\end{proof}

Let now $f\in C_{b}^{4}(S^{1}\times S^{1},\widetilde{\mathbb{R}})^{H}$ be a
cocycle satisfying $c([f])=0$. Define $h_{1}:(S^{1}\times S^{1})^{4}%
\rightarrow\mathbb{R}$ as%
\[
h_{1}((x_{1},y_{1}),...,(x_{4},y_{4}))=\left\{
\begin{array}
[c]{l}%
f((x_{i},y_{j}),(x_{1},y_{1}),...,(x_{4},y_{4})),\\
\text{ \ \ \ \ \ \ \ if }\exists\text{ }i\neq i^{\prime},\text{ }j\neq
j^{\prime}\text{ with }x_{i}=x_{i^{\prime}}\text{, }y_{j}=y_{j^{\prime}},\\
0,\text{ \ \ \ \ otherwise.}%
\end{array}
\right.
\]
Let us check that $h_{1}$ is well defined: If the condition $x_{i}%
=x_{i^{\prime}}$ is satisfied for different pairs $i_{1}\neq i_{1}^{\prime}$
and $i_{2}\neq i_{2}^{\prime}$ (thus $\{i_{1},i_{1}^{\prime}\}\neq
\{i_{2},i_{2}^{\prime}\}$), then $\sharp\{x_{1},...,x_{4}\}\leq2$ and
\[
n_{1}((x_{i_{\ell}},y_{j}),(x_{1},y_{1}),...,(x_{4},y_{4}))\leq2,\text{ \ for
}\ell=1,2.
\]
By Step \ref{Step 0}, this now implies that $f$ vanishes on both of those
$5$-tuples and hence $h_{1}((x_{1},y_{1}),...,(x_{4},y_{4}))=0$ is well
defined. The case when the condition $y_{j}=y_{j^{\prime}}$ is satisfied for
different pairs of indices is treated symmetrically. Observe furthermore that
$h_{1}$ belongs to $C_{b}^{3}(S^{1}\times S^{1},\widetilde{\mathbb{R}})^{H}$
because $f$ belongs to $C_{b}^{4}(S^{1}\times S^{1},\widetilde{\mathbb{R}%
})^{H}$.

\begin{step}
\label{Step 1}Set $f_{1}=f-\delta h_{1}\in C_{b}^{4}(S^{1}\times
S^{1},\widetilde{\mathbb{R}})^{H}$. If $\underline{z}\in(S^{1}\times
S^{1})^{5}$ satisfies $n_{1}(\underline{z})=n_{2}(\underline{z})=3$, then
$f_{1}(\underline{z})=0$.
\end{step}

\begin{proof}
Because $f$ is alternating and $H$-invariant, and since $\mathrm{Isom}%
(\mathbb{H}^{2})$ acts transitively on oriented triples of distinct points of
$S^{1}$, the value of $f$ on $5$-tuples $\underline{z}$ with $n_{1}%
(\underline{z})=n_{2}(\underline{z})=3$ only depends on the configuration of
the coordinates of $\underline{z}$. There are, up to permutation, five such
configurations. Thus, there exists $\lambda_{0},...,\lambda_{4}$ in
$\mathbb{R}$ such that for every triple $(x_{0},x_{1},x_{2})$ and
$(y_{0},y_{1},y_{2})$ of distinct points of $S^{1}$, the following equalities
hold:%
\begin{align*}
f((x_{0},y_{0}),(x_{0},y_{1}),(x_{1},y_{0}),(x_{1},y_{1}),(x_{2},y_{2}))  &
=\lambda_{0}\cdot\mathrm{Or}(x_{0},x_{1},x_{2})\cdot\mathrm{Or}(y_{0}%
,y_{1},y_{2}),\\
f((x_{0},y_{0}),(x_{1},y_{0}),(x_{2},y_{0}),(x_{2},y_{1}),(x_{2},y_{2}))  &
=\lambda_{1}\cdot\mathrm{Or}(x_{0},x_{1},x_{2})\cdot\mathrm{Or}(y_{0}%
,y_{1},y_{2}),\\
-f((x_{0},y_{0}),(x_{1},y_{0}),(x_{1},y_{1}),(x_{2},y_{1}),(x_{2},y_{2}))  &
=\lambda_{2}\cdot\mathrm{Or}(x_{0},x_{1},x_{2})\cdot\mathrm{Or}(y_{0}%
,y_{1},y_{2}),\\
f((x_{0},y_{0}),(x_{1},y_{0}),(x_{1},y_{1}),(x_{1},y_{2}),(x_{2},y_{2}))  &
=\lambda_{3}\cdot\mathrm{Or}(x_{0},x_{1},x_{2})\cdot\mathrm{Or}(y_{0}%
,y_{1},y_{2}),\\
f((x_{0},y_{0}),(x_{0},y_{1}),(x_{1},y_{1}),(x_{2},y_{1}),(x_{2},y_{2}))  &
=\lambda_{4}\cdot\mathrm{Or}(x_{0},x_{1},x_{2})\cdot\mathrm{Or}(y_{0}%
,y_{1},y_{2}).
\end{align*}
Note that by Step \ref{Step 0}, the above relations also hold when
$n_{1}(\underline{z})\leq2$ or $n_{2}(\underline{z})\leq2$ since both sides of
the equations are then equal to $0$.

We start by invoking Lemma \ref{Lemma: c(z)=signs c(z)} to show that
$\lambda_{0}=0$: The even permutation exchanging the first with the second and
the third with the fourth coordinate of
\[
((x_{0},y_{0}),(x_{0},y_{1}),(x_{1},y_{0}),(x_{1},y_{1}),(x_{2},y_{2}))
\]
is realized on the first factor by the identity and on the second by the
reversing orientation isometry permuting $y_{0}$ with $y_{1}$ and fixing
$y_{2}$. In particular, $f$ has to vanish on this $5$-tuple.

Furthermore, note that $\lambda_{3}=\lambda_{4}$ since $f$ is invariant under
the orientation preserving isometry $\tau$ of $\mathbb{H}^{2}\times
\mathbb{H}^{2}$ permuting the two factors.

From the cocycle relation $\delta f=0$, we compute%
\begin{align*}
0  &  =\delta f((x_{0},y_{0}),(x_{1},y_{0}),(x_{1},y_{1}),(x_{2},y_{1}%
),(x_{2},y_{2}),(x_{2},y_{0}))\\
&  =(\lambda_{1}+\lambda_{2}-\lambda_{3}-\lambda_{4})\mathrm{Or}(x_{0}%
,x_{1},x_{2})\cdot\mathrm{Or}(y_{0},y_{1},y_{2}),
\end{align*}
and we see that%
\begin{equation}
\lambda_{1}+\lambda_{2}=\lambda_{3}+\lambda_{4}=2\lambda_{3}.
\label{Equ: lambda1+lambda2=2xlambda3}%
\end{equation}

\begin{claim}
If $f\in C_{b}^{4}(S^{1}\times S^{1},\widetilde{\mathbb{R}})^{H\text{ }}$ is
such that $c([f])=0$, then%
\[
2\left(  \lambda_{1}+\lambda_{2}\right)  +\lambda_{3}+\lambda_{4}=2\left(
\lambda_{1}+\lambda_{2}+\lambda_{3}\right)  =0.
\]

\end{claim}

\begin{proof}
[Proof of Claim]Let%
\[
\Gamma_{2}=\left\langle \left.  a_{1},b_{1},a_{2},b_{2}\right\vert \text{
}[a_{1},b_{1}][a_{2},b_{2}]=1\right\rangle <\mathrm{PSL}_{2}\mathbb{R}%
\]
be a representation of the fundamental group $\Gamma_{2}$ of the genus $2$
surface $\Sigma_{2}$ in $\mathrm{PSL}_{2}\mathbb{R}$. Note that%
\begin{align*}
z  &  =(1,a_{1},b_{1})+(1,a_{1}b_{1},a_{2})+(1,a_{1}b_{1}a_{2},b_{2})\\
&  \text{ \ \ \ \ }-(1,b_{2},a_{2})-(1,b_{2}a_{2},b_{1})-(1,b_{2}a_{2}%
b_{1},a_{1})
\end{align*}
is a cycle in $C_{2}(\Gamma_{2})\hookrightarrow C_{2}(\mathrm{PSL}%
_{2}\mathbb{R)}$ representing the fundamental class $[\Sigma_{2}]\in
H_{2}(\Sigma_{2})\cong H_{2}(\Gamma_{2})$. Recall that given two $2$-chains
$(g_{0},g_{1},g_{2})$ and $(k_{0},k_{1},k_{2})$ in $C_{2}(\mathrm{PSL}%
_{2}\mathbb{R)}$, their product $(g_{0},g_{1},g_{2})\times(k_{0},k_{1},k_{2})$
in $C_{4}(\mathrm{PSL}_{2}\mathbb{R\times}\mathrm{PSL}_{2}\mathbb{R)}$ is
defined as the $4$-chain%
\begin{align*}
&  \text{ \ \ }((g_{0},k_{0}),(g_{0},k_{1}),(g_{0},k_{2}),(g_{1},k_{2}%
),(g_{2},k_{2}))\\
&  -((g_{0},k_{0}),(g_{0},k_{1}),(g_{1},k_{1}),(g_{1},k_{2}),(g_{2},k_{2}))\\
&  +((g_{0},k_{0}),(g_{0},k_{1}),(g_{1},k_{1}),(g_{2},k_{1}),(g_{2},k_{2}))\\
&  +((g_{0},k_{0}),(g_{1},k_{0}),(g_{1},k_{1}),(g_{1},k_{2}),(g_{2},k_{2}))\\
&  -((g_{0},k_{0}),(g_{1},k_{0}),(g_{1},k_{1}),(g_{2},k_{1}),(g_{2},k_{2}))\\
&  +((g_{0},k_{0}),(g_{1},k_{0}),(g_{2},k_{0}),(g_{2},k_{1}),(g_{2},k_{2})).
\end{align*}
Thus, $z\times z$ is a $4$-cycle in $C_{4}(G)\hookrightarrow C_{4}(H)$.

For any cocycle $f$ in $C_{b}^{4}(S^{1}\times S^{1},\widetilde{\mathbb{R}%
})^{H\text{ }}$, the cohomology class $c([f])\in H_{cb}^{4}(H,\widetilde
{\mathbb{R}})$ is represented in $C^{4}(H,\widetilde{\mathbb{R}})$ by the
cocycle
\[%
\begin{array}
[c]{rccl}%
f_{\xi}: & H^{5} & \longrightarrow & \mathbb{R}\\
& (h_{0},...,h_{4}) & \longmapsto & f(h_{0}(\xi,\xi),...,h_{4}(\xi,\xi)),
\end{array}
\]
where $\xi$ is a fixed base point in $S^{1}$. For any $(g_{0},g_{1},g_{2})$
and $(k_{0},k_{1},k_{2})$ in $C_{2}(\mathrm{PSL}_{2}\mathbb{R)}$, we have%
\begin{equation}
f_{\xi}((g_{0},g_{1},g_{2})\times(k_{0},k_{1},k_{2}))=(2(\lambda_{1}%
+\lambda_{2})+\lambda_{3}+\lambda_{4})\mathrm{Or}(g_{0}\xi,g_{1}\xi,g_{2}%
\xi)\cdot\mathrm{Or}(k_{0}\xi,k_{1}\xi,k_{2}\xi).
\label{Equ: f_xi(product of (gi's x ki's))}%
\end{equation}

Upon conjugating $\Gamma_{2}$, we can without loss of generality assume that
$a_{1}\xi=\xi$. Now, remember that, as seen in Section
\ref{Section: volume form}, $\mathrm{Or}_{\xi}$ is a cocycle in $C^{2}%
(\mathrm{PSL}_{2}\mathbb{R},\mathbb{R})$ representing $(1/\pi)\omega
_{\mathbb{H}^{2}}$. In particular, its evaluation on the fundamental class
$[\Sigma_{2}]$ is equal to $(1/\pi)\cdot\mathrm{Vol}(\Sigma_{2})=4$, so that%
\begin{align}
4  &  =\left\langle \mathrm{Or}_{\xi},z\right\rangle =\mathrm{Or}(\xi,a_{1}%
\xi,b_{1}\xi)+\mathrm{Or}(\xi,a_{1}b_{1}\xi,a_{2}\xi)+\mathrm{Or}(\xi
,a_{1}b_{1}a_{2}\xi,b_{2}\xi)\nonumber\\
&  \text{ \ \ \ \ \ \ \ \ \ \ \ \ \ \ \ \ \ }-\mathrm{Or}(\xi,b_{2}\xi
,a_{2}\xi)-\mathrm{Or}(\xi,b_{2}a_{2}\xi,b_{1}\xi)-\mathrm{Or}(\xi,b_{2}%
a_{2}b_{1}\xi,a_{1}\xi)\nonumber\\
&  =\mathrm{Or}(\xi,a_{1}b_{1}\xi,a_{2}\xi)+\mathrm{Or}(\xi,a_{1}b_{1}a_{2}%
\xi,b_{2}\xi)-\mathrm{Or}(\xi,b_{2}\xi,a_{2}\xi)-\mathrm{Or}(\xi,b_{2}a_{2}%
\xi,b_{1}\xi), \label{Equ: (Or,z)= 4 times orientation cocycle}%
\end{align}
since $a_{1}\xi=\xi$. Because the cocycle $\mathrm{Or}$ takes its values in
$\{-1,0,+1\}$, it is now immediate that%
\[
\mathrm{Or}(\xi,a_{1}b_{1}\xi,a_{2}\xi)=\mathrm{Or}(\xi,a_{1}b_{1}a_{2}%
\xi,b_{2}\xi)=-\mathrm{Or}(\xi,b_{2}\xi,a_{2}\xi)=-\mathrm{Or}(\xi,b_{2}%
a_{2}\xi,b_{1}\xi)=1.
\]
Note that alternatively, the above equalities can be checked directly by
studying the action of $\Gamma_{2}$ on $\partial\mathbb{H}^{2}$.

Finally, the assumption that $c([f])=0$ tells us that $f_{\xi}$ is a
coboundary and hence vanishes on cycles. In particular, we get $f_{\xi
}(z\times z)=0$. But from $\left(  \ref{Equ: f_xi(product of (gi's x ki's))}%
\right)  $ and $\left(  \ref{Equ: (Or,z)= 4 times orientation cocycle}\right)
$, we straightforwardly compute%
\[
f_{\xi}(z\times z)=16\cdot\left(  2(\lambda_{1}+\lambda_{2})+\lambda
_{3}+\lambda_{4}\right)  ,
\]
which proves the claim.
\end{proof}

Denote by $\lambda_{0}^{\prime},...,\lambda_{4}^{\prime}$ the real numbers in
the defining equations for $f$ on the $5$-tuples $\underline{z}$ with
$n_{1}(\underline{z})=n_{2}(\underline{z})=3$ we would obtain by replacing $f$
by $f_{1}$. Note that, as for $f$, we have $\lambda_{0}^{\prime}=0$. From the
definition of $f_{1}$ as $f-\delta h_{1}$, we furthermore obtain%
\[%
\begin{array}
[c]{ll}%
\lambda_{1}^{\prime}=\lambda_{1}-\lambda_{1}=0, & \text{ \ \ }\lambda
_{2}^{\prime}=\lambda_{1}+\lambda_{2},\\
\lambda_{3}^{\prime}=\lambda_{3}, & \text{ \ \ }\lambda_{4}^{\prime}%
=\lambda_{4}.
\end{array}
\]
But from $\left(  \ref{Equ: lambda1+lambda2=2xlambda3}\right)  $ and the
claim, it now follows that $\lambda_{j}^{\prime}=0$, for every $0\leq j\leq4$,
which proves that $f_{1}$ vanishes on all $5$-tuples $\underline{z}$ with
$n_{1}(\underline{z})=n_{2}(\underline{z})=3$.
\end{proof}

Define $h_{2}:(S^{1}\times S^{1})^{4}\rightarrow\mathbb{R}$ as%
\[%
\begin{array}
[c]{l}%
h_{2}((x_{1},y_{1}),...,(x_{4},y_{4}))=\\
\text{ \ \ }=\left\{
\begin{array}
[c]{l}%
\frac{1}{2}\left[  f_{1}((x_{i},y_{k}),(x_{1},y_{1}),...,(x_{4},y_{4}%
))+f_{1}((x_{i},y_{\ell}),(x_{1},y_{1}),...,(x_{4},y_{4}))\right]  ,\\
\text{ \ \ \ \ \ \ if }\{i,j,k,\ell\}=\{1,2,3,4\}\text{ and }x_{i}=x_{j},\\
\frac{1}{2}\left[  f_{1}((x_{k},y_{i}),(x_{1},y_{1}),...,(x_{4},y_{4}%
))+f_{1}((x_{\ell},y_{i}),(x_{1},y_{1}),...,(x_{4},y_{4}))\right]  ,\\
\text{ \ \ \ \ \ \ if }\{i,j,k,\ell\}=\{1,2,3,4\}\text{ and }y_{i}=y_{j},\\
0,\text{ \ \ \ otherwise.}%
\end{array}
\right.
\end{array}
\]
To check that $h_{2}$ is well defined, we verify that if the first or the
second condition are verified by different sets of indices, then $h_{2}$ is in
both cases defined as $0$: If the condition $x_{i}=x_{j}$ is satisfied for
different pairs $i_{1}\neq j_{1}$ and $i_{2}\neq j_{2}$ (thus $\{i_{1}%
,j_{1}\}\neq\{i_{2},j_{2}\}$), then as in the proof that $h_{1}$ is well
defined, we get%
\[
n_{1}((x_{i_{\ell}},y_{\ast}),(x_{1},y_{1}),...,(x_{4},y_{4}))\leq2,
\]
for $\ell=1,2$ and $y_{\ast}\in\{y_{1},...,y_{4}\}$. By Step \ref{Step 0},
this implies that $f_{1}$ evaluated on those $5$-tuples vanishes, and
$h_{2}((x_{1},y_{1}),...,(x_{4},y_{4}))=0$ is well defined. The case when the
condition $y_{i}=y_{j}$ is satisfied for different pairs of indices is treated
symmetrically. Finally, suppose that $\{i,j,k,\ell\}=\{i^{\prime},j^{\prime
},k^{\prime},\ell^{\prime}\}=\{1,2,3,4\}$ with $x_{i}=x_{j}$ and
$y_{i^{\prime}}=y_{j^{\prime}}$. Then we have both%
\[
\sharp\{x_{1},...,x_{4}\}\leq3\text{ \ and \ }\sharp\{y_{1},...,y_{4}\}\leq3.
\]
In particular, both $n_{1}$ and $n_{2}$ are at most equal to $3$ when
evaluated on the $5$-tuples appearing in the definition of $h_{2}$. Since by
Step \ref{Step 1}, $f_{1}$ vanishes on those $5$-tuples, we obtain, in this
case also, that $h_{2}((x_{1},y_{1}),...,(x_{4},y_{4}))=0$ is well defined.
Observe furthermore that $h_{2}$ belongs to $C_{b}^{3}(S^{1}\times
S^{1},\widetilde{\mathbb{R}})^{H}$ because $f_{1}$ belongs to $C_{b}^{4}%
(S^{1}\times S^{1},\widetilde{\mathbb{R}})^{H}$.

\begin{step}
\label{Step 2}Set $f_{2}=f_{1}-\delta h_{2}\in C_{b}^{4}(S^{1}\times
S^{1},\widetilde{\mathbb{R}})^{H}$. If $\underline{z}\in(S^{1}\times
S^{1})^{5}$ satisfies $n_{1}(\underline{z})+n_{2}(\underline{z})\leq7$, then
$f_{2}(\underline{z})=0$.
\end{step}

\begin{proof}
By Step \ref{Step 0}, $f_{2}(\underline{z})=0$ whenever $n_{1}(\underline
{z})\leq2$ or $n_{2}(\underline{z})\leq2$. If $n_{1}(\underline{z}%
)=n_{2}(\underline{z})=3$, then $f_{1}(\underline{z})=0$ by Step \ref{Step 1}.
Furthermore, in this case $\delta h_{2}(\underline{z})$ is also equal to $0$
since all the $5$-tuples $\underline{z}^{\prime}$ evaluated on by $f_{1}$ in
the definition of $h_{2}(\underline{\widehat{z}^{i}})$, where $\underline
{\widehat{z}^{i}}$ denotes the $4$-tuple obtained from $\underline{z}$ by
removing its $i$-th coordinate, for $0\leq i\leq4$, satisfy $n_{1}%
(\underline{z}^{\prime})=n_{2}(\underline{z}^{\prime})=3$. Thus, $f_{2}$
vanishes in this case also. By symmetry, it now remains to treat the case
$n_{1}(\underline{z})=3$, $n_{2}(\underline{z})=4$.

Up to permutation, we have two possibilities for the first factor:

\begin{enumerate}
\item $x_{0}=x_{1}=x_{2}\neq x_{3}\neq x_{4},$

\item $x_{0}=x_{1}\neq x_{2}=x_{3}\neq x_{4}$.
\end{enumerate}

\begin{enumerate}
\item In the first case, we then have, again up to permutation, three options
for the second factor:

\begin{enumerate}
\item $y_{0}=y_{1}$: Trivially, $f_{2}(\underline{z})=0$ since the two first
coordinates of $\underline{z}$ are equal.

\item $y_{0}=y_{4}$: We consider two subcases:

\begin{itemize}
\item $\left\langle y_{0},y_{3}\right\rangle \cap\left\langle y_{1}%
,y_{2}\right\rangle \neq\varnothing$: By Lemma
\ref{Lemma: 4 points, geodesics intersect or not}, there exists an orientation
preserving isometry $h$ of $\mathbb{H}^{2}$ exchanging $y_{0}$ with $y_{3}$
and $y_{1}$ with $y_{2}$. Furthermore, there exists an orientation reversing
isometry $g$ of $\mathbb{H}^{2}$ with $gx_{0}=x_{0}$, $gx_{3}=x_{4}$ and
$gx_{4}=x_{3}$. Since $f_{2}$ is alternating and $H$-invariant, we get on the
one hand, applying the even permutation $(1$ $2)(3$ $4)$ and the action by
$(g,h)$,
\begin{align}
&  \text{\ }f_{2}((x_{0},y_{0}),(x_{0},y_{1}),(x_{0},y_{2}),(x_{3}%
,y_{3}),(x_{4},y_{0}))\nonumber\\
\text{ \ \ \ \ \ \ \ \ \ \ \ \ \ \ \ \ \ \ \ \ \ \ }  &  =f_{2}((x_{0}%
,y_{0}),(x_{0},y_{2}),(x_{0},y_{1}),(x_{4},y_{0}),(x_{3},y_{3}))\nonumber\\
&  =-f_{2}((gx_{0},hy_{0}),(gx_{0},hy_{2}),(gx_{0},hy_{1}),(gx_{4}%
,hy_{0}),(gx_{3},hy_{3}))\nonumber\\
&  =-f_{2}((x_{0},y_{3}),(x_{0},y_{1}),(x_{0},y_{2}),(x_{3},y_{3}%
),(x_{4},y_{0})). \label{Equ in Step 2 case 1.2.1}%
\end{align}
On the other hand, apply the cocycle relation of $f_{2}$ to the $6$-tuple%
\[
\text{ \ \ \ \ \ \ \ \ \ \ \ \ \ \ \ \ \ \ \ \ \ \ }\underline{w}%
=((x_{0},y_{3}),(x_{0},y_{0}),(x_{0},y_{1}),(x_{0},y_{2}),(x_{3},y_{3}%
),(x_{4},y_{0})).
\]
If one removes the $3$-rd or the $4$-th variable of $\underline{w}$, then the
remaining $5$-tuple has $n_{1}=n_{2}=3$ and thus $f_{2}$ vanishes on it. If
one removes the $5$-th or the $6$-th variable of $\underline{w}$, then the
remaining $5$-tuple has $n_{1}=2$ and here also $f_{2}$ vanishes on it. The
cocycle relation $\delta f_{2}(\underline{w})$ hence simplifies to
\begin{align*}
\text{\ \ \ }  &  \text{ }f_{2}((x_{0},y_{0}),(x_{0},y_{1}),(x_{0}%
,y_{2}),(x_{3},y_{3}),(x_{4},y_{0}))\\
\text{\ \ \ \ \ \ \ \ \ }\text{ \ \ \ \ \ \ \ \ \ \ \ \ }  &  =f_{2}%
((x_{0},y_{3}),(x_{0},y_{1}),(x_{0},y_{2}),(x_{3},y_{3}),(x_{4},y_{0})).
\end{align*}
Together with $\left(  \ref{Equ in Step 2 case 1.2.1}\right)  $, this shows
that $f_{2}$ vanishes on $5$-tuples of the form $((x_{0},y_{0}),(x_{0}%
,y_{1}),(x_{0},y_{2}),(x_{3},y_{3}),(x_{4},y_{0})).$

\item $\left\langle y_{0},y_{3}\right\rangle \cap\left\langle y_{1}%
,y_{2}\right\rangle =\varnothing$: By Lemma
\ref{Lemma: 4 points, geodesics intersect or not}, there exists an orientation
reversing isometry $h$ of $\mathbb{H}^{2}$ exchanging $y_{0}$ with $y_{3}$ and
$y_{1}$ with $y_{2}$. As above, there exists an orientation reversing isometry
$g$ of $\mathbb{H}^{2}$ with $gx_{0}=x_{0}$, $gx_{3}=x_{4}$ and $gx_{4}=x_{3}%
$. Since $f_{1}$ is alternating and $H$-invariant, we get, applying the even
permutation $(1$ $2)(3$ $4)$ and the action by $(g,h)$,
\begin{align*}
\text{ \ \ \ \ \ \ }  &  \text{ }f_{1}((x_{0},y_{0}),(x_{0},y_{1}%
),(x_{0},y_{2}),(x_{3},y_{3}),(x_{4},y_{0}))\\
\text{\ \ \ \ \ \ \ \ \ }\text{ \ \ \ \ \ \ \ \ \ \ \ \ }  &  =f_{1}%
((x_{0},y_{0}),(x_{0},y_{2}),(x_{0},y_{1}),(x_{4},y_{0}),(x_{3},y_{3}))\\
\text{\ \ \ \ \ \ \ \ \ }\text{ \ \ \ \ \ }  &  =f_{1}((gx_{0},hy_{0}%
),(gx_{0},hy_{2}),(gx_{0},hy_{1}),(gx_{4},hy_{0}),(gx_{3},hy_{3}))\\
\text{\ \ \ \ \ \ \ \ \ }\text{ \ \ \ \ \ }  &  =f_{1}((x_{0},y_{3}%
),(x_{0},y_{1}),(x_{0},y_{2}),(x_{3},y_{3}),(x_{4},y_{0})).
\end{align*}
In particular,%
\begin{align*}
\text{ \ \ \ \ \ \ \ \ \ \ \ }\text{ }  &  \delta h_{2}((x_{0},y_{0}%
),(x_{0},y_{1}),(x_{0},y_{2}),(x_{3},y_{3}),(x_{4},y_{0}))\\
\text{ }  &  =h_{2}((x_{0},y_{1}),(x_{0},y_{2}),(x_{3},y_{3}),(x_{4},y_{0}))\\
&  =\frac{1}{2}\left[  f_{1}((x_{0},y_{3}),(x_{0},y_{1}),(x_{0},y_{2}%
),(x_{3},y_{3}),(x_{4},y_{0}))\right. \\
&  \text{ \ \ }\left.  +f_{1}((x_{0},y_{0}),(x_{0},y_{1}),(x_{0},y_{2}%
),(x_{3},y_{3}),(x_{4},y_{0}))\right] \\
&  =f_{1}((x_{0},y_{0}),(x_{0},y_{1}),(x_{0},y_{2}),(x_{3},y_{3}),(x_{4}%
,y_{0})),
\end{align*}
and hence $f_{2}$ vanishes on this $5$-tuple.
\end{itemize}

\item $y_{3}=y_{4}$: Our $5$-tuple $\underline{z}$ has the form
\[
\underline{z}=((x_{0},y_{0}),(x_{0},y_{1}),(x_{0},y_{2}),(x_{3},y_{3}%
),(x_{4},y_{3})).
\]
Set
\[
\underline{w}=(\underline{z},(x_{0},y_{3})).
\]
The cocycle relation $\delta f_{2}(\underline{w})=0$ gives
\[
f_{2}(\underline{z})=\sum_{i=0}^{4}(-1)^{i}f_{2}(\underline{\widehat{w}}%
^{i}),
\]
where $\underline{\widehat{w}}^{i}$ denotes the $5$-tuple obtained from
$\underline{w}$ by removing its $i$-th coordinate, for $0\leq i\leq4$. But
$f_{2}(\underline{\widehat{w}}^{i})=0$ for $i=0,1,2$, since in this case
$n_{1}(\underline{\widehat{w}}^{i})=n_{2}(\underline{\widehat{w}}^{i})=3$, and
for $i=3,4,$ since then $n_{1}(\underline{\widehat{w}}^{i})=2$. In particular,
$f_{2}(\underline{z})=0$.
\end{enumerate}

\item In the second case, we have up to permutation, three options for the
second factor:

\begin{enumerate}
\item $y_{0}=y_{1}$:\textbf{ }Again, trivially, $f_{2}(\underline{z})=0$.

\item $y_{0}=y_{4}$:\textbf{ }Our $5$-tuple $\underline{z}$ has the form
\[
\underline{z}=((x_{0},y_{0}),(x_{0},y_{1}),(x_{2},y_{2}),(x_{2},y_{3}%
),(x_{4},y_{0})).
\]
Set
\[
\underline{w}=(\underline{z},(x_{2},y_{0})).
\]
The cocycle relation $\delta f_{2}(\underline{w})=0$ gives
\[
f_{2}(\underline{z})=\sum_{i=0}^{4}(-1)^{i}f_{2}(\underline{\widehat{w}}%
^{i}).
\]
But $f_{2}(\underline{\widehat{w}}^{i})=0$ for $i=1,2,3$, since in this case
$n_{1}(\underline{\widehat{w}}^{i})=n_{2}(\underline{\widehat{w}}^{i})=3$, and
for $i=4,$ since then $n_{1}(\underline{\widehat{w}}^{i})=2$. Finally, for
$i=0$, we have $(n_{1}(\underline{\widehat{w}}^{i}),n_{2}(\underline
{\widehat{w}}^{i}))=(3,4)$, but $\underline{\widehat{w}}^{i}$ is of the form
treated in $(1)$ since its first coordinates consists of the $5$-tuple
$(x_{0},x_{2},x_{2},x_{4},x_{2})$ and hence vanishes when evaluated on $f_{2}%
$. It follows that $f_{2}(\underline{z})=0$.

\item $y_{0}=y_{2}$: Our $5$-tuple $\underline{z}$ has the form
\[
\underline{z}=((x_{0},y_{0}),(x_{0},y_{1}),(x_{2},y_{0}),(x_{2},y_{3}%
),(x_{4},y_{4})).
\]
Set
\[
\underline{w}=(\underline{z},(x_{2},y_{4})).
\]
The cocycle relation $\delta f_{2}(\underline{w})=0$ gives
\[
f_{2}(\underline{z})=\sum_{i=0}^{4}(-1)^{i}f_{2}(\underline{\widehat{w}}%
^{i}).
\]
We see that: for $i=1,3,$ $n_{1}(\underline{\widehat{w}}^{i})=n_{2}%
(\underline{\widehat{w}}^{i})=3$; for $i=4$, $n_{1}(\underline{\widehat{w}%
}^{i})=2$; for $i=0$, $(n_{1}(\underline{\widehat{w}}^{i}),n_{2}%
(\underline{\widehat{w}}^{i}))=(3,4)$, but $\underline{\widehat{w}}^{i}$ is of
the form treated in $(1)$ since its first coordinates consists of the
$5$-tuple $(x_{0},x_{2},x_{2},x_{4},x_{2})$; for $i=2$, again $(n_{1}%
(\underline{\widehat{w}}^{i}),n_{2}(\underline{\widehat{w}}^{i}))=(3,4)$, but
\[
\underline{\widehat{w}}^{0}=((x_{0},y_{0}),(x_{0},y_{1}),(x_{2},y_{3}%
),(x_{4},y_{4}),(x_{2},y_{4}))
\]
is of the form treated in $(2b)$. In all those cases, we thus obtain
$f_{2}(\underline{\widehat{w}}^{i})=0$ and hence $f_{2}(\underline{z})=0$.
\end{enumerate}
\end{enumerate}
\end{proof}

Define $h_{3}:(S^{1}\times S^{1})^{4}\rightarrow\mathbb{R}$ as%
\[
h_{3}((x_{1},y_{1}),...,(x_{4},y_{4}))=\frac{1}{12}\sum_{i,j=1}^{4}%
f_{2}((x_{i},y_{j}),(x_{1},y_{1}),...,(x_{4},y_{4})).
\]
Observe furthermore that $h_{3}$ belongs to $C_{b}^{3}(S^{1}\times
S^{1},\widetilde{\mathbb{R}})^{H}$ because $f_{2}$ belongs to $C_{b}^{4}%
(S^{1}\times S^{1},\widetilde{\mathbb{R}})^{H}$.

\begin{step}
\label{Step 3}Set $f_{3}=f_{2}-\delta h_{3}\in C_{b}^{4}(S^{1}\times
S^{1},\widetilde{\mathbb{R}})^{H}$. For any $\underline{z}\in(S^{1}\times
S^{1})^{5}$, we have $f_{3}(\underline{z})=0$.
\end{step}

\begin{proof}
We start with a preliminary computation.

\begin{claim}
For any $\left(  (x_{1},y_{1}),...,(x_{4},y_{4})\right)  $ in $\left(
S^{1}\times S^{1}\right)  ^{4}$ and any $1\leq i\neq j\leq4$ and $1\leq
i^{\prime}\neq j^{\prime}\leq4$, we have%
\[
f_{2}((x_{i},y_{j}),(x_{1},y_{1}),...,(x_{4},y_{4}))=f_{2}((x_{i^{\prime}%
},y_{j^{\prime}}),(x_{1},y_{1}),...,(x_{4},y_{4})).
\]

\end{claim}

\begin{proof}
[Proof of Claim]By symmetry, if we prove the claim for $i=i^{\prime}$, then it
is also proven for $j=j^{\prime}$. Furthermore, the general case then follows
since if $i\neq j^{\prime}$, we can go from $(i,j)$ to $(i,j^{\prime})$ and
then to $(i^{\prime},j^{\prime})$, and similarly if $i^{\prime}\neq j$.
Finally, if $i=j^{\prime}$ and $i^{\prime}=j$, there exists $k\neq i,j$ so
that we can go from $(i,j)$ to $(i,k)$ to $(i^{\prime},k)$ and to $(i^{\prime
},j^{\prime})$. Thus, it is now enough to prove the claim for $i=i^{\prime}$.

Consider the $6$-tuple
\[
\underline{z}=((x_{i},y_{j}),(x_{i^{\prime}},y_{j^{\prime}}),(x_{1}%
,y_{1}),...,(x_{4},y_{4})).
\]
From the cocycle relation $\delta f_{2}(\underline{z})$ we see that the claim
would follow from the equality%
\begin{equation}
0=\sum_{k=1}^{4}(-1)^{k}f_{2}((x_{i},y_{j}),(x_{i^{\prime}},y_{j^{\prime}%
}),(x_{1},y_{1}),...,\widehat{(x_{k},y_{k})},...,(x_{4},y_{4})).
\label{Equation: in claim in step 3}%
\end{equation}
But $n_{1}$ evaluated on those $5$-tuple is smaller or equal to $4$ when $k=i$
and to $3$ otherwise, while $n_{2}$ is smaller or equal to $4$ when $k=j$ or
$j^{\prime}$ and to $3$ otherwise. Because $k$ can not simultaneously be equal
to $i$ and $j$ or $j^{\prime}$ it follows that $n_{1}+n_{2}\leq7$ on all of
the $5$-tuples appearing in $\left(  \ref{Equation: in claim in step 3}%
\right)  $, so that, by Step \ref{Step 2}, each of the summand in $\left(
\ref{Equation: in claim in step 3}\right)  $ is equal to $0$.
\end{proof}

Note that it follows that
\begin{equation}
h_{3}((x_{1},y_{1}),...,(x_{4},y_{4}))=f_{2}((x_{i},y_{j}),(x_{1}%
,y_{1}),...,(x_{4},y_{4})), \label{Equ: h3= f2 for any i,j}%
\end{equation}
for any $1\leq i\neq j\leq4$.

Observe that by Step \ref{Step 2}, both $f_{2}$ and $\delta h_{3}$ vanish on
$5$-tuples $\underline{z}$ satisfying $n_{1}(\underline{z})+n_{2}%
(\underline{z})\leq7$, so that the same holds for $f_{3}$. We now will prove
step by step, that $f_{3}$ also vanishes on $5$-tuples $\underline{z}$ with
$(n_{1}(\underline{z}),n_{2}(\underline{z}))=(3,5),(4,4),(4,5)$ and $(5,5)$.
In all but one subcase, the strategy is the same as in most of the proof of
Step \ref{Step 2}: 1) Start with an arbitrary $5$-tuple with given
$(n_{1}(\underline{z}),n_{2}(\underline{z}))$. 2) Apply the cocycle relation
$\delta f_{3}=0$ to an appropriately chosen $6$-tuple $\underline
{w}=(\underline{z},(x,y))$, for $(x,y)$ in $S^{1}\times S^{1}$, so that
\[
f_{3}(\underline{z})=\sum_{i=0}^{4}(-1)^{i}f_{3}(\underline{\widehat{w}}%
^{i}),
\]
where $\underline{\widehat{w}}^{i}$ denotes, for $0\leq i\leq5$, the $5$-tuple
obtained from $\underline{w}$ by removing its $i$-th coordinate. (In
particular, $\underline{\widehat{w}}^{5}=\underline{z}$.) 3) Show that
$f_{3}(\underline{\widehat{w}}^{i})=0$ for $0\leq i\leq4$.

\begin{description}
\item[(3,5)] Let $\underline{z}$ be a $5$-tuple with $n_{1}(\underline{z})=3$
and $n_{2}(\underline{z})=5$. We distinguish two subcases:

\begin{itemize}
\item The $5$-tuple $\underline{z}$ has, up to permutation, the form
\[
\underline{z}=((x_{0},y_{0}),(x_{1},y_{1}),(x_{2},y_{2}),(x_{2},y_{3}%
),(x_{2},y_{4})).
\]
Set
\[
\underline{w}=(\underline{z},(x_{2},y_{0})).
\]
We have $f_{3}(\underline{\widehat{w}}^{i})=0$ for $i=0,1$ because
$n_{1}(\underline{\widehat{w}}^{i})=2$, and for $i=2,3,4$ because
$n_{1}(\underline{\widehat{w}}^{i})+n_{2}(\underline{\widehat{w}}^{i})=3+4$.

\item The $5$-tuple $\underline{z}$ has, up to permutation, the form
\[
\underline{z}=((x_{0},y_{0}),(x_{0},y_{1}),(x_{2},y_{2}),(x_{2},y_{3}%
),(x_{4},y_{4})).
\]
Set
\[
\underline{w}=(\underline{z},(x_{0},y_{4})).
\]
We have $n_{1}(\underline{\widehat{w}}^{i})=3$ for $i=0,1,2,3$ and
$n_{1}(\underline{\widehat{w}}^{i})=2$ for $i=4$. Also, $n_{2}(\underline
{\widehat{w}}^{i})=4$ for $i=0,1,2,3$ and $n_{2}(\underline{\widehat{w}}%
^{i})=5$ for $i=4$. In any case, $n_{1}(\underline{\widehat{w}}^{i}%
)+n_{2}(\underline{\widehat{w}}^{i})=7$ so that $f_{3}(\underline{\widehat{w}%
}^{i})=0$ for every $i$.
\end{itemize}

Note that it follows that $f_{3}(\underline{z})=0$ whenever $n_{1}%
(\underline{z})=3$ or $n_{2}(\underline{z})=3$.

\item[(4,4)] We distinguish two subcases:

\begin{itemize}
\item The $5$-tuple $\underline{z}$ has, up to permutation, the form%
\[
\underline{z}=((x_{1},y_{2}),(x_{1},y_{1}),...,(x_{4},y_{4})).
\]
In this case we have
\[
\delta h_{3}(\underline{z})=h_{3}((x_{1},y_{1}),...,(x_{4},y_{4}%
))=f_{2}(\underline{z}),
\]
where the last equality follows from $\left(  \ref{Equ: h3= f2 for any i,j}%
\right)  $. In particular, $f_{3}(\underline{z})=f_{2}(\underline{z})-\delta
h_{3}(\underline{z})=0$, as desired.

\item The $5$-tuple $\underline{z}$ has, up to permutation, the form%
\[
\underline{z}=((x_{0},y_{0}),(x_{0},y_{1}),(x_{2},y_{2}),(x_{3},y_{2}%
),(x_{4},y_{4})).
\]
Set
\[
\underline{w}=(\underline{z},(x_{0},y_{2})).
\]
We have $f_{3}(\underline{\widehat{w}}^{i})=0$ for $i=0,1$ because
$n_{2}(\underline{\widehat{w}}^{i})=3$, and for $i=2,3,4$ because
$n_{1}(\underline{\widehat{w}}^{i})=3$.
\end{itemize}

\item[(4,5)] We can assume that $\underline{z}$ has the form%
\[
\underline{z}=((x_{0},y_{0}),(x_{0},y_{1}),(x_{2},y_{2}),(x_{3},y_{3}%
),(x_{4},y_{4})).
\]
Set
\[
\underline{w}=(\underline{z},(x_{3},y_{4})).
\]
We have $f_{3}(\underline{\widehat{w}}^{i})=0$ for $i=0,1,3$ because
$n_{1}(\underline{\widehat{w}}^{i})=n_{2}(\underline{\widehat{w}}^{i})=4$, and
for $i=2,4$ because $n_{1}(\underline{\widehat{w}}^{i})=3$.

\item[(5,5)] Let finally $\underline{z}$ be a generic $5$-tuple. Set
\[
\underline{w}=(\underline{z},(x_{3},y_{4})).
\]
We have $f_{3}(\underline{\widehat{w}}^{i})=0$ for $i=0,1,2,3$ because
$n_{2}(\underline{\widehat{w}}^{i})=4$, and for $i=4$ because $n_{1}%
(\underline{\widehat{w}}^{i})=4$.
\end{description}
\end{proof}

\textbf{In conclusion, }the arbitrary cocycle $f\in C_{b}^{4}(S^{1}\times
S^{1},\widetilde{\mathbb{R}})^{H}$ satisfying $c([f])=0$ we started with is a
coboundary since
\[
f=f_{1}+\delta h_{1}=f_{2}+\delta h_{2}+\delta h_{1}=f_{3}+\delta h_{3}+\delta
h_{2}+\delta h_{1}=\delta(h_{1}+h_{2}+h_{3}),
\]
and Theorem \ref{Thm: comparison map} is hence proven.

\end{document}